\documentclass[11pt,a4paper]{scrartcl}

\usepackage{amsmath,amsfonts, amssymb,  mathrsfs}
\usepackage{array}			

\usepackage{tikz}
\usetikzlibrary{shapes,arrows,shadows}

\newcommand{\R}{{\mathbb{R}}}       
\newcommand{\E}{\operatorname{E}}
\newcommand{\Cov}{\operatorname{Cov}}
\newcommand{\MSE}{\operatorname{MSE}}

\newcommand{\T}{^{\top}}

\newcommand{\st}{\operatorname{s.t.}}

\newcommand{\LL}{L^{2}(\Omega)}

\usepackage{setspace}

\title{Efficient response surface methods based on generic surrogate models\thanks{This research has been supported by the German Federal Ministry of Economics and Technology (BMWi) within the joint project ComFliTe.}}
\author{
Benjamin Rosenbaum\footnotemark[2]\and Volker Schulz\footnotemark[2]
}
\date{}

\begin{document}

\maketitle

\renewcommand{\thefootnote}{\fnsymbol{footnote}}
\footnotetext[2]{Universit\"at Trier, FB IV -
Abteilung Mathematik, Universit\"atsring 15, D-54296 Trier, Germany (benjamin.rosenbaum@uni-trier.de, volker.schulz@uni-trier.de)}
\renewcommand{\thefootnote}{\arabic{footnote}}


%
%
%


\onehalfspacing
\parskip3ex

\textbf{Abstract.} Surrogate models are used for global approximation of responses generated by expensive computer experiments like CFD applications. In this paper, we make use of structural similarities which are shared by a class of related problems. We identify these structures by applying statistical shape models. They are used to build a generic surrogate model approximation to sample data of a new problem of the same class. In a variable fidelity framework the generic surrogate model is combined with the sample data to generate an efficient and globally accurate interpolation model, which requires less costly sample evaluations than ordinary response surface methods. We demonstrate our method with an aerodynamic test case and show that it significantly improves the approximation quality.

\section{Introduction}\label{sec:intro} 

In multidisciplinary numerical simulation and optimization, surrogate modeling has gained popularity during the last two decades. Numerical computation of realistic models still is challenging and computationally intensive. When the global behavior of an input-output relationship of a computer experiment is sought, dense evaluations over the whole  input parameter space are out of reach. Surrogate models or also called metamodels or response surface models approximate or interpolate the output of a computer code based on a moderate number of evaluations. Typically, the evaluation of the computer experiment dominates the overall computational cost of model generation. For the sake of efficiency, methods have to be designed to use as few samples as possible without sacrificing accuracy. Radial basis functions and the Kriging method are widely used because of their flexibility and ability to interpolate highly nonlinear functions. Improving these methods is a major topic of ongoing research: On the one hand, adaptive sampling strategies are investigated to reduce the number of required evaluations \cite{ShanW10},\cite{BR2011}. On the other hand, additional information which is assumed to be correlated with the response is used to improve the accuracy of approximation, like gradient enhanced Kriging (GEK) \cite{LaurenceauS08},\cite{YamazakiRM2010},\cite{LockwoodA12}, Cokriging \cite{ForresterK09},\cite{HanZG2010},\cite{HanZG12} and variable fidelity modeling (VFM) \cite{HanGH2010},\cite{MooreRP11}.

All these approaches treat the response itself more or less as an unknown output of a black-box. This paper is motivated by the assumption that for a predefined problem class, the behavior of the response is not arbitrary, but rather related to other instances of the mutual problem class. For example in computational fluid dynamics (CFD), responses of aerodynamic coefficients, depending on the input parameters Mach number and angle of attack, share structural similarities for different airfoil geometries. To identify these structures we make use of the concept of statistical shape models, which use a principal component analysis to quantify modes of variation of a previously computed training database. If functions of a problem class are accessible in form of such a database, a new test case of this function family can be approximated using only few evaluations. Instead of directly interpolating the samples, first the principal components are fitted to the data which then act as a generic surrogate model (GSM). Based on this model, interpolation of the samples is performed in a VFM framework. 

The statistical interpolation method Kriging has been widely used in geostatistics since the 1950s \cite{Krige1951},\cite{Matheron1963}, it was introduced in surrogate modeling for computer experiments in the eighties \cite{SacksWMW89} and is nowadays applied to CFD simulation and optimization, see e.g.\ \cite{LaurenceauS08},\cite{LamKHP09},\cite{ForresterSK08}. A survey about sampling strategies in optimization is given in \cite{ShanW10}, while adaptive sampling strategies for global approximation can be found in \cite{BusbyFI07},\cite{LaurenceauS08},\cite{CrombecqTGD09},\cite{GramacyL09}. The authors in \cite{JinCS02} present an early review and a recent one is given in \cite{BR2011}.  Another framework of improving the approximation quality uses secondary information about the response. Gradient enhanced Kriging (GEK) incorporates derivatives \cite{MorrisMY93},\cite{KoehlerO96}, which are often available by adjoint computations. This approach is used in various fields of aerodynamic applications such as optimization \cite{ChungA02},\cite{Liu03},\cite{LaurenceauM08},\cite{YamazakiRM2010}, uncertainty quantification \cite{DwightH09},\cite{LockwoodA12} and global approximation \cite{LaurenceauS08},\cite{BR2011}. GEK can be interpreted as a special form of Cokriging. This method also originates from geostatistics \cite{Wackernagel03} and enables the incorporation of any auxiliary variable which is correlated with the primary one. In the last few years, it was also applied in aerodynamics and engineering \cite{ForresterSK2007},\cite{ForresterK09},\cite{HanZG2010}, where the approximation of few costly computations of high fidelity could be improved by a large number of evaluations of a cheaper low fidelity model. Other variable fidelity modeling (VFM) methods use bridge functions to correct the discrepancy between data of low and high fidelity and date back almost as long as surrogate modeling for computer experiments itself, e.g.\ \cite{ChangHGK93}. Nowadays they are also used in aerodynamics \cite{TangGL05},\cite{RobinsonEWH06},\cite{HanGH2010},\cite{KozielY11}. Recently a new VFM method was developed \cite{HanGoertz2011}, which overcomes the difficulties of model building and robustness in Cokriging as well as the problems of accuracy and missing mean squared error prediction in the correction based VFM methods. While statistical shape models are popular in the fields of computer vision and (medical) image processing \cite{Cootesetal1994},\cite{LameckerSHD04},\cite{DaviesTT2008}, to our knowledge this paper represents the first attempt to use this technique in response surface methods for global approximation of expensive computer experiments. 

The outline is as follows. Section \ref{sec:rsm} presents Kriging as a surrogate modeling framework for computer experiments. A database of surrogates for a predefined problem class from aerodynamic simulation is considered  in section \ref{sec:database} and the problem of establishing a correspondence between the database elements is discussed. The concept of statistical shape models as well as the POD method for $L^{2}$-functions are presented in section \ref{sec:POD}. In section \ref{sec:GPOD}, we describe the gappy POD method, extend it to the continuous case and explain how the gappy POD fit of the POD basis elements to some sample data establishes a generic surrogate model. We discuss how the generic surrogate model approximation and the sample data are combined by a variable fidelity modeling interpolation in section \ref{sec:HK}. In the final section, we present numerical results of the generic surrogate modeling technique for responses of aerodynamic coefficients with a database of several airfoil geometries and compare them to common Kriging.     

\section{Surrogate modeling}\label{sec:rsm} 

We consider a computer experiment like a CFD solver, whose evaluation is expensive. It is treated as a black-box, depending on input parameters $x\in\R^{d}$ and we observe a scalar \textit{response} $y(x)$. We want to approximate the unknown function $y:\R^d\rightarrow\R$ in a domain of interest $\Omega\subset\R^{d}$ by a surrogate model $\widehat y(x)$ which can be evaluated at low computational cost. The Kriging method is a statistical method of interpolating $n$ given data pairs $\left\{(x_{i},y(x_{i}))\right\}_{i=1}^{n}$ of evaluations of the unknown function $y(x)$. In this section we only give a brief review, following \cite{SacksWMW89} and \cite{KoehlerO96}. For more information about constructing the surrogate we refer to \cite{DACE02},\cite{ForresterSK08},\cite{SantnerWN03}. For the interpolation model, the deterministic response is treated as a sum of a linear regression part and a ``lack of fit'' term
\begin{equation}
  y(x)=\sum_{k=1}^{K}\beta_{k}f_{k}(x)+z(x), \quad x\in\Omega\subset\R^{d}, \label{sec1:KrigingModel}
\end{equation}
where $f_{k}(x):\R^{d}\rightarrow\R$ are known functions with coefficients $\beta_{k}\in\R$ and $z(x)$ is the realization of a stationary Gaussian process, which captures the nonlinear behavior of the response.  The regression part usually consists of low order polynomials (simplest case: $K$=1, $f_{1}\equiv 1$) and the Gaussian process is assumed to satisfy 
\begin{equation}
  \E\left[z(x)\right]=0, \quad \Cov\left[z(x),z(w)\right]=\sigma^{2}R(w-x) \quad \forall x,w\in\Omega, \label{sec1:process}
\end{equation}
where $\sigma^{2}$ is the process variance and $R(w-x)$ is a spatial correlation function, which will be explained further at the end of this section.

The Kriging predictor 
\begin{equation}
  \widehat y(x):=c(x)^{\top}Y \label{sec1:Kriging}
\end{equation}
is a weighted sum of the given response evaluations $Y:=\left(y(x_{1}),\dots,y(x_{n})\right)^{\top}$. With the weights $c_{i}(x)$ chosen as the solution of the optimization problem 
\begin{equation}
\begin{split}
  \min_{c(x)\in\R^{n}}\ &\MSE[\widehat{y}(x)]=\E[(c(x)^{\top}Y-y(x))^2]\\
  \operatorname{s.t.}\ & \E[c(x)^{\top}Y-y(x)]=0, \label{sec1:BLUE}
\end{split}
\end{equation}
it is also a \textit{best linear unbiased estimator} (BLUE). The solution of (\ref{sec1:BLUE}) is given by the solution of the linear equation 
\begin{equation}
  \begin{bmatrix}R & F \\F\T & 0 \end{bmatrix} \begin{pmatrix} c(x) \\ \mu(x) \end{pmatrix} = \begin{pmatrix} r(x) \\ f(x)\end{pmatrix} \label{sec1:krigingsystem}
\end{equation}
where $R=\left[R(x_i,x_j)\right]_{i,j}\in\R^{n\times n}$ denotes the positive definite and symmetric correlation matrix, $r(x)=\left(R(x,x_i)\right)_i\in\R^n$ contains the correlations between $x$ and every $x_i$. $F=\left[f_k(x_i)\right]_{i,k}\in\R^{n\times K}$ is the linear regression design matrix and $f(x)=\left(f_k(x)\right)_k\in\R^K$ contains the evaluations of the design functions in $x$. $\mu(x)\in\R^{K}$ are the Lagrange multipliers for the unbiasedness condition in (\ref{sec1:BLUE}).

Using (\ref{sec1:krigingsystem}), a closed term for the Kriging estimator (\ref{sec1:Kriging}) can now easily be derived:
\begin{equation}
  \widehat y(x)=c(x)^{\top}Y = \begin{pmatrix} r(x) \\ f(x)\end{pmatrix}^{\top} \begin{bmatrix}R & F \\F\T & 0 \end{bmatrix}^{-1}\begin{pmatrix} Y \\ 0\end{pmatrix}. \label{sec1:closedform}
\end{equation}
Inserting any $x_{i}$ $(i=1,\dots,n)$, it becomes obvious that the Kriging predictor is in fact an interpolator ($\widehat y(x_{i})=y(x_{i})$). 
Solving the linear equation independently from $x$ in (\ref{sec1:closedform}) only once for a given dataset, the Kriging predictor can be evaluated efficiently at the cost of a dot product of size $n+K$, ending up with a cheaply accessible surrogate model for $y(x)$.

The spatial correlation function (\ref{sec1:process}) is modeled as a product of one dimensional correlation functions
\begin{equation}
  R(w-x,\theta)=\prod_{k=1}^{d}R^{(k)}(\bigl|w^{(k)}-x^{(k)}\bigr|,\theta^{(k)}), 
\end{equation} 
still depending on so-called hyperparameters $\theta=\left(\theta^{(1)},\dots,\theta^{(d)} \right)$ which determine the correlation lengths. Popular choices for the correlation function are exponential functions of the type ${R^{(k)}(\bigl|w^{(k)}-x^{(k)}\bigr|,\theta^{(k)})}=\exp\left\{-\theta^{(k)}\bigl|w^{(k)}-x^{(k)}\bigr|^{p}\right\}$, $p\in\left[1,2\right]$, or cubic splines, all satisfying $R(0)=1$ and $R(\left|w-x\right|)\xrightarrow[\left|w-x\right|\rightarrow\infty]{} 0$. The hyperparameters $\theta^{(k)}$ are  usually determined by solving a \textit{maximum likelihood problem}
\begin{equation}
  \max_{\theta,\beta,\sigma^{2}}\ (2\pi)^{-\frac{n}{2}}\sigma^{-n}(\det R(\theta))^{-\frac{1}{2}}\exp \left\{-\frac{1}{2\sigma^2}(Y-F\beta)\T R(\theta)^{-1}(Y-F\beta) \right\} .\label{sec1:MLE1}
\end{equation}
$R(\theta)$ again denotes the correlation matrix, now depending on $\theta$. It is possible to cancel out $\beta$ and $\sigma^{2}$ by necessary first order conditions 
\begin{align}
  \beta(\theta) &= \left( F\T R(\theta)^{-1} F\right)^{-1}F\T R(\theta)^{-1}Y \label{sec1:MLbeta} \\
  \sigma^2\left(\theta,\beta(\theta)\right) &=\frac{1}{n}(Y-F\beta(\theta))\T R(\theta)^{-1}(Y-F\beta(\theta)) \label{sec1:MLsigma}
\end{align}
and the problem can be reduced to 
\begin{equation}
  \min_\theta\ \left\{ \sigma^2\left(\theta,\beta\left(\theta\right)\right) \left( \det R(\theta) \right)^\frac{1}{n} \right\}. \label{sec1:MLE2}
\end{equation}
However, when solving the maximum likelihood problem with a gradient based algorithm, formulation (\ref{sec1:MLE2}) can cause difficulties due to highly nonlinear implicit dependencies and it is advisable to rather use (\ref{sec1:MLE1}) \cite{BR2011}.

\section{Function database for aerodynamic simulation}\label{sec:database} 

We now concretize our test case for better understanding without loss of generality. Our goal is to generate surrogate models for scalar aerodynamic coefficients lift ($c_{l}$), drag ($c_{d}$) and pitching moment ($c_{m}$) for two-dimensional airfoil geometries, which depend on the input parameters Mach number ($\textit{Ma}$) and angle of attack ($\alpha$). For every input parameter configuration, the response can be evaluated by running a simulation with a (computationally intensive) CFD solver. We assume that the responses of different airfoils constitute a mutual problem class and share structural similarities. When already having generated surrogate models for several airfoil geometries, theses similarities can be used to improve the quality of a surrogate model for a new instance of this problem class. 

We consider a database of surrogate models $y_{1}(x),\dots,y_{m}(x)$ for the same problem class and we omit the superscript for simplicity ($y_{i}:=\widehat y_{i}$). In our case, each $y_{i}\in C^{1}(\Omega)$ is the response of an aerodynamic coefficient (e.g.\ lift) depending on the input parameters $x=(\textit{Ma},\alpha)$ for a particular airfoil geometry. So the database consists of previously computed response surface functions corresponding to $m$ different airfoils. Furthermore we assume each response function to have a sufficient accuracy, i.e.\ it represents a globally valid surrogate model in the domain of interest $\Omega\subset\R^{2}$.

In the fields of computer vision and image processing, \textit{statistical shape models} built from a dataset of examples have been widely used \cite{DaviesTT2008}. After establishing correspondence between the database functions by admissible transformations, also called alignment, a principal component analysis is performed to identify the most important modes of variation. Often Euclidean or similarity transformations are used to establish correspondence between the images,  e.g.\ see \cite{Cootesetal1994}. For an overview on image registration techniques and possible transformations see \cite{FischerM2008},\cite{Szeliski2006},\cite{ZitovaF2003}.


As an \textit{admissible transformation} we define the following one, which is an affine transformation in each dimension of $\Omega\subset\R^{2}$ and in the image space $y(\Omega)\subset\R$:
\begin{align}
  \overline x(q^{j})&:=\begin{pmatrix} x^{(1)}(1+q^{j}_{1})+q^{j}_{2} \\ x^{(2)}(1+q^{j}_{3})+q^{j}_{4}\end{pmatrix}, \label{secdatabase:xtransform} \\
  \overline y_{j}\left(\overline x(q^{j}),q^{j} \right) &:= y_{j}\left(\overline x(q^{j})\right)(1+q^{j}_{5})+q^{j}_{6}, \label{secdatabase:ytransform}
\end{align}
parametrized by $q^{j}\in\R^{6}$ ($j=1,\dots,m$). Note that the transformation is chosen such that $\overline y_{j}\left(\overline x(0),0 \right)=y_{j}(x)$ for $q^{j}=0$. We emphasize that this transformation was found to be suitable for our test cases, other problem classes could require other admissible transformations. One function $y_{1}(x)$ is defined as a reference, meaning no transformation is applied ($q^{1}:=0$). For the other transformation parameters $q^{2},\dots,q^{m}$ an optimization problem with $6(m-1)$ unknowns has to be solved:
\begin{equation}
\min_{q^{2},\dots,q^{m}} \frac{1}{m(m-1)}\sum_{j=1}^{m}\sum_{k>j}\int_{\Omega}\left( \overline y_{j}\bigl(\overline x(q^{j}),q^{j} \bigr) - \overline y_{k}\bigl(\overline x(q^{k}),q^{k} \bigl)  \right)^{2} dx + \frac{\delta}{2}\sum_{j=2}^{m}{q^{j}}^{\top}q^{j}. \label{secdatabase:optproblem}
\end{equation}
The solution minimizes the overall sum of squared differences, while for robustness a penalty term is included which guarantees that the transformation does not become too large. This nonlinear least squares problem is solved by a Gau{\ss}-Newton algorithm \cite{NocedalW06}. 

Some computational issues will now be addressed briefly. The integral in (\ref{secdatabase:optproblem}) must be approximated by a numerical quadrature
\begin{equation}
  \int_{\Omega}f(x)dx\approx \sum_{i=1}^{N}w_{i}f(x_{i}),
\end{equation}
e.g.\ with the $x_{i}$ as elements of a $\sqrt{N}\times\sqrt{N}$-tensorgrid and positive weights $w_{i}$.   With the term
\begin{gather}
  e_{i,j,k}:=\sqrt{w_{i}}\left( \overline y_{j}\bigl(\overline x(x_{i},q^{j}),q^{j} \bigr) - \overline y_{k}\bigl(\overline x(x_i,q^{k}),q^{k} \bigl)  \right),\\
  e\in\R^{N\frac{m(m-1)}{2}}, 
\end{gather}
we can write the first summand of (\ref{secdatabase:optproblem}) in discretized form as a sum of squared differences
\begin{equation}
  SSD(q)=\frac{1}{m(m-1)}\sum_{k=1}^{m}\sum_{j>k}\sum_{i=1}^{N}e_{i,j,k}^{2}.
\end{equation}
Then using
\begin{equation}
  J:=\frac{\partial e(q)}{\partial q}\in\R^{N\frac{m(m-1)}{2}\times 6m},  
\end{equation}
the Gau{\ss}-Newton algorithm needs the gradient and an approximation to the Hessian
\begin{align}
  \nabla_{q}SSD(q)&=J^{\top}e  \in\R^{6m}\\
  \nabla^{2}_{q}SSD(q)&\approx J^{\top}J  \in\R^{6m\times 6m}.  
\end{align}
Each entry of the symmetric matrix $J^{\top}J$ is a dot product of length $N\frac{m(m-1)}{2}$ and there are $\frac{6m(6m+1)}{2}$ entries to compute, so the algorithm can really benefit from parallelization. Furthermore, the memory usage by storing the matrix $J$ is $\mathcal O(Nm^{3})$.  

\section{Proper orthogonal decomposition}\label{sec:POD} 

For the dataset $\left(y_{1},\dots,y_{m}\right)$, $y_{i}\in L^{2}(\Omega;\R)$, we want to perform a proper orthogonal decomposition (POD) which is also called principal component analysis (PCA) or Karhunen-Lo\`{e}ve transformation, see also \cite{BerkoozHL93} or \cite{Sirovich87}. Again, we omit the superscripts for simplicity ($y_{i}:=\overline{\widehat y}_{i}$, see (\ref{sec1:Kriging}),(\ref{secdatabase:ytransform})). We keep in mind that in this paper, all functions are surrogate models which have been aligned by admissible transformations, but the POD method applies for any $y_{i}\in\LL$ as well. In PCA significant structures of the dataset are identified which is realized by an orthogonal decomposition of the covariance matrix 
\begin{equation}
  C^{(y)}:=\left[ \left(y_{i},y_{j} \right)_{\LL}\right]_{i,j}\in \R^{m\times m}. \label{secPOD:covariancematrix}
\end{equation}
We point out that classically mean centered functions are considered, i.e.\ the proper orthogonal decomposition is performed on $\left(\check y_{1},\dots,\check y_{m}\right)$, $\check y_{i}:=y_{i}-\frac{1}{m}\sum_{k=1}^{m}y_{k}$ \cite{DaviesTT2008}. This leads to a decomposition of the space of variations from the mean instead of the space spanned by the functions themselves. In the statistical shape model
\begin{equation}
  \frac{1}{m}\sum_{k=1}^{m}y_{k}+\sum_{j=1}^{l}a_{j}\psi_{j}
\end{equation}
the principal components of variation $\psi_{j}\in\LL$ from the mean are controlled by parameters $a_{j}\in\R$ $(j=1,\dots,l)$. We investigated proper orthogonal decompositions of both mean centered $\left(\check y_{1},\dots,\check y_{m}\right)$ and plain datasets $\left(y_{1},\dots,y_{m}\right)$. In our test cases, mean centering did not produce better results while increasing the complexity of algorithmic implementation, so we decided to use the plain formulation.

In this section, we introduce the POD method for the continuous case and follow the 
discussion of \cite{Volkwein01}. Let $\Omega\subset\R^{d}$ be bounded, $L^{2}(\Omega;\R)$ denotes the Hilbert space with $\left(f,g\right)_{\LL}=\int_{\Omega}f(x)g(x)dx$ and $\|f\|_{\LL}^{2}=\int_{\Omega}f(x)^{2}dx$ for all $f,g\in\LL$. We assume the $m$ functions $y_{i}\in\LL$ $(i=1,\dots,m)$ to be linear independent and define $y(x):=\left(y_{1}(x),\dots,y_{m}(x)\right)\in\R^{m}$. For the POD method, consider $\mathcal Y^{m}=\operatorname{span}\left\{y_{1},\dots,y_{m}\right\}\subset\LL$ of dimension $m$. Let $\left\{\psi_{1},\dots,\psi_{m}\right\}$ be an orthonormal basis of $\mathcal Y^{m}$. Clearly, when using the complete basis, every $y_{i}$ can be represented by a linear combination of its elements
\begin{equation}
  y_{i}=\sum_{j=1}^{m}\left(y_{i},\psi_{j}\right)_{\LL}\psi_{j} \quad(i=1,\dots,m).
\end{equation}
But we search for an orthonormal set of functions $\left\{\psi_{1},\dots,\psi_{l}\right\}\subset \mathcal Y^{m}$ of dimension $l\leq m$, which describes $\mathcal Y^{m}$ as good as possible, i.e.\ every $y_{i}$ is approximated by a linear combination of $\left\{\psi_{1},\dots,\psi_{l}\right\}$. This leads to the following optimization problem:
\begin{align}
\begin{split}
  \min_{\psi_{1},\dots,\psi_{l}\in \mathcal Y^{m}} & \sum_{i=1}^{m} \biggl\|y_{i}-\sum_{j=1}^{l}\left(y_{i},\psi_{j} \right)_{\LL}\psi_{j} \biggr\|_{\LL}^{2} \\
  \st \quad & \left(\psi_{i},\psi_{j}\right)_{\LL}=\delta_{ij} \quad (i,j=1,\dots,l)
\end{split} \label{secPOD:prob1}
\end{align} 
Expanding the norm in (\ref{secPOD:prob1}) and using $\left(\psi_{i},\psi_{j}\right)_{\LL}=\delta_{ij}$ yields an equivalent formulation:
\begin{align}
\begin{split}
  \max_{\psi_{1},\dots,\psi_{l}\in \mathcal Y^{m}} & \sum_{i=1}^{m} \sum_{j=1}^{l}\left(y_{i},\psi_{j} \right)_{\LL}^{2} \\
  \st \quad & \left(\psi_{i},\psi_{j}\right)_{\LL}=\delta_{ij} \quad (i,j=1,\dots,l)
\end{split}\label{secPOD:prob2}
\end{align} 
The solution $\left\{\psi_{1},\dots,\psi_{l}\right\}$ of (\ref{secPOD:prob2}) is called POD-basis of rank $l$. We now briefly explain how the POD-basis is determined.

We define the operator $\mathcal C:\LL\rightarrow \mathcal Y^{m}$,
\begin{equation}
  \mathcal C \psi := \sum_{i=1}^{m}\left(\psi,y_{i} \right)_{\LL}y_{i}.
\end{equation}
Then for $\mathcal C$ exists a series of orthonormal eigenfunctions $\left\{\psi_{i} \right\}_{i=1}^{\infty}$ and corresponding nonnegative real eigenvalues $\left\{\lambda_{i} \right\}_{i=1}^{\infty}$ with
\begin{gather}
  \mathcal C \psi_{i}=\lambda_{i}\psi_{i}  \label{secPOD:eigensystem} \\
  \lambda_{1}\geq \ldots \geq \lambda_{m} > 0 \\
  \lambda_{i} = 0 \quad (i > m),
\end{gather}
see \cite{KunischV2003}. Furthermore, $\left\{\psi_{i} \right\}_{i=1}^{\infty}$ is an orthonormal basis for $\LL$ and $\operatorname{span}\{ \psi_{1},\ldots,\psi_{m}\}=\mathcal Y^{m}$, which implicates that $\{ \psi_{1},\dots,\psi_{m}\}$ is an orthonormal basis for $\mathcal Y^{m}$. Also, for any $l\leq m$ $\{ \psi_{1},\dots,\psi_{l}\}$ is the unique solution of (\ref{secPOD:prob1}). The error can be expressed by the sum of the $m-l$ remaining eigenvalues 
\begin{equation}
  \sum_{i=1}^{m} \biggl\|y_{i}-\sum_{j=1}^{l}\left(y_{i},\psi_{j} \right)_{\LL}\psi_{j} \biggr\|_{\LL}^{2} = \sum_{j=l+1}^{m}\lambda_{j}. \label{secPOD:errorformula}
\end{equation}

So the POD-basis is determined by the first $l$ eigenfunctions of $\mathcal C$ which correspond to the $l$ largest eigenvalues. A finite approach for solving the eigensystem (\ref{secPOD:eigensystem}) can be derived \cite{Volkwein01}: Because $\psi_{k}\in \mathcal Y^{m}$ $(k=1,\dots,l)$, we set
\begin{equation}
  \psi_{k}=\kappa_{k}\sum_{i=1}^{m}v_{i}^{k}y_{i}\quad (k=1,\dots,l). \label{secPOD:lin.comb}
\end{equation}
Inserting (\ref{secPOD:lin.comb}) into the eigensystem (\ref{secPOD:eigensystem}) yields
\begin{equation}
  \kappa_{k}\sum_{j=1}^{m}\left( \sum_{i=1}^{m}\left( y_{i},y_{j} \right)_{\LL} v_{i}^{k}\right)y_{j} = \kappa_{k} \sum_{j=1}^{m}\left( \lambda_{k}v_{j}^{k} \right)y_{j} \quad (k=1,\dots,l)
\end{equation}
and exploiting the linear independence of $\left\{y_{1},\dots,y_{m}\right\}$ we conclude
\begin{equation}
  \sum_{i=1}^{m}\left( y_{i},y_{j} \right)_{\LL} v_{i}^{k}=  \lambda_{k}v_{j}^{k}  \quad (k=1,\dots,l;\ j=1,\dots,m ).
\end{equation}
This is a discrete eigenvalue problem and with the covariance matrix $C=C^{(y)}$ from (\ref{secPOD:covariancematrix}) we can write it as
\begin{equation}
  Cv^{k}=\lambda_{k}v^{k} \quad (k=1,\dots,l).
\end{equation}
So for determining the POD-basis $\{ \psi_{1},\dots,\psi_{l}\}$ of rank $l$ one has to compute the $l$ eigenvectors $v^{k}\in\R^{m}$ of $C$ which correspond to the $l$ largest eigenvalues $\lambda_{k}$. Setting $\kappa_{k}=\frac{1}{\sqrt{\lambda_{k}}}$ for normalization in (\ref{secPOD:lin.comb}), the POD-basis elements are
\begin{equation}
  \psi_{k}=\frac{1}{\sqrt{\lambda_{k}}}\sum_{i=1}^{m}v_{i}^{k}y_{i}\quad (k=1,\dots,l). 
\end{equation}
With notations $V_{l}:=\left[v_{i}^{k} \right]_{i,k}\in\R^{m\times l}$ and $\Sigma_{l}:=\operatorname{diag}(\sqrt{\lambda_{1}},\dots,\sqrt{\lambda_{l}})\in\R^{l\times l}$ we can also write
\begin{equation}
  \psi(x):=(\psi_{1}(x),\dots,\psi_{l}(x))=y(x)V_{l}\Sigma_{l}^{-1}\in \R^{1\times l}. \label{secPOD:closedterm}
\end{equation}


\section{Gappy POD in Hilbert spaces}\label{sec:GPOD} 

The gappy POD method was first introduced in \cite{EversonS1995} for reconstructing images of faces from incomplete data. A mask function was used which set the greyscale of every pixel to zero which was not part of the incomplete dataset. In \cite{BuiDW2004}, the gappy POD method was applied to CFD problems for the first time and a selection vector was used to set the missing flow solution vector's entries to zero. We now introduce a straightforward approach for applying this methodology to $\LL$, where the gappy data is given by function evaluations on a discrete subset of $\Omega$. But first we briefly recap how a ``non gappy'' function is approximated by the POD-basis.

When approximating an arbitrary function $\phi\in\LL$ (not necessarily $\phi\in\mathcal{Y}^m$) with the POD-basis, the solution of the optimization problem 
\begin{equation}
  \min_{a_{1}^{(\psi)},\dots,a_{l}^{(\psi)}\in\R}\frac{1}{2}\biggl\| \phi(x)-\sum_{j=1}^{l}a_j^{(\psi)}\psi_{j}(x) \biggr\|_{\LL}^{2} \label{secGPOD:minL2}
\end{equation}
is sought. This is a linear least squares problem and using optimality conditions and exploiting $\left( \psi_{i},\psi_{j} \right)_{\LL}=\delta_{ij}$ one derives
\begin{equation}
  a_{j}^{(\psi)}=\left( \phi,\psi_{j} \right)_{\LL} \quad (j=1,\dots,l).
\end{equation}
The approximating function $\widetilde\phi(x):=\sum_{j=1}^{l}a_j^{(\psi)}\psi_{j}(x)$ is a projection of $\phi(x)$ onto the subspace $\mathcal Y^{l}:=\operatorname{span}\left\{\psi_{1},\dots,\psi_{l}\right\}\subset\mathcal Y^{m}\subset\LL$.

Now for the gappy POD method, suppose $\phi(x)$ itself is unknown again and only a set of $n$ data pairs of evaluations
\begin{gather}
  \left\{ (x_{i};\phi(x_{i})) \right\}_{i=1}^{n},\\
  x_{i}\in\Omega\subset\R^{d},\ \phi(x_{i})\in\R,\ x_{i}\neq x_{j}\ (i\neq j),\quad (i=1,\dots,n) \nonumber
\end{gather}
is available. Based on the data, we want to reconstruct the unknown function $\phi\in\LL$. Assuming there is a strong relation between $\phi$ and $\mathcal Y^{m}$, we approximate $\phi$ by a $\widetilde\phi\in\mathcal Y^{l}$. Similar to (\ref{secGPOD:minL2}), we pose an optimization problem
\begin{equation}
  \min_{a_{1}^{(\psi)},\dots,a_{l}^{(\psi)}\in\R}\frac{1}{2}\sum_{i=1}^{n} \biggl( \phi(x_{i})-\sum_{j=1}^{l}a_j^{(\psi)}\psi_{j}(x_{i}) \biggr)^{2} \label{secGPOD:minsum},
\end{equation}
which again is a linear least squares problem. Clearly, though the $\psi_{j}$ are orthonormal in $\LL$, meaning $(\psi_{i},\psi_{j})_{\LL}=\delta_{ij}$, this is not the case anymore on a subset of $\Omega$. Particularly, $\sum_{k=1}^{n}\psi_{i}(x_{k})\psi_{j}(x_{k})$ $(i\neq j)$ is generally not equal to zero. Introducing the design matrix
\begin{equation}
  \Psi:=\left[ \psi_{j}\left(x_{i}\right) \right]_{i=1,j=1}^{n,l}\in\R^{n\times l},\ n\geq l,\ \operatorname{rank}\Psi =l.
\end{equation}
and $\varphi:=\left(\phi(x_{1}),\dots,\phi(x_{n})\right)\T\in\R^{n\times 1}$, the solution of (\ref{secGPOD:minsum}) is given by the linear equation
\begin{equation}
  \Psi^{\top}\Psi a^{(\psi)}= \Psi^{\top}\varphi. \label{secGPOD:linsystem}
\end{equation}

Since every $\psi_{j}$ $(j=1,\dots,l)$ is a linear combination of $\left\{ y_{1},\dots,y_{m}\right\}$, we now derive a framework for avoiding multiple redundant evaluations of $y_{i}(x)$ both while solving (\ref{secGPOD:linsystem}) and also for evaluating $\widetilde\phi(x)=\sum_{j=1}^{l}a_{j}^{(\psi)}\psi_{j}(x)$. Even if each $y_{j}(x)$ is a surrogate model, depending on its complexity and $l$, $m$ and $n$, the evaluations needed can be a bottleneck in the model generation.  We define 
\begin{equation}
  \bold{Y} := \left[ y_{j}\left(x_{i}\right) \right]_{i=1,j=1}^{n,m}\in\R^{n\times m}
\end{equation}
and analogously to (\ref{secPOD:closedterm}) we get
\begin{equation}
  \Psi = \bold{Y}V_{l}\Sigma_{l}^{-{1}}.
\end{equation}
Also, with (\ref{secPOD:closedterm}) we get for the approximating function 
\begin{align}
  \widetilde\phi(x)&=\sum_{j=1}^{l}a_{j}^{(\psi)}\psi_{j}(x) \nonumber \\ 
  &= \psi(x)a^{(\psi)} \nonumber \\
  &= y(x)\underbrace{V_{l}\Sigma_{l}^{-1}a^{(\psi)}}_{:=a^{(y)}\in\R^{m\times 1}} \label{secGPOD:approxfun}
\end{align} 
such that we have a closed form of $\widetilde \phi (x)=y(x)a^{(y)}=\sum_{j=1}^{m}a_{j}^{(y)}y_{j}(x)$ as a linear combination of $\left\{y_{1},\dots,y_{m} \right\}$.

Because the database elements $y_{1},\dots,y_{m}$ have been transformed by solving the alignment problem (\ref{secdatabase:optproblem}), the unknown function $\phi(x)$ (respectively its known evaluations $\left\{ (x_{i};\phi(x_{i})) \right\}_{i=1}^{n}$) has to be allowed a transformation of the same class. At this point the task is not finding an overall alignment, but fitting the POD basis to a fixed set of data pairs $\left\{ (x_{i},\phi(x_{i})) \right\}_{i=1}^{n}$. So rather than transforming $\phi(x)$, we (equivalently) apply a transformation to the approximating function $\widetilde\phi(x)$, parametrized by $p$. The double bar indicates that we deal with a second transformation after the already known initial transformation of the database functions by alignment (parametrized by $q$):
\begin{align}
  \overline{\overline x}(p)&:=\begin{pmatrix}\overline x^{(1)}(q^{j})(1+p_{1})+p_{2}\\ \overline x^{(2)}(q^{j})(1+p_{3})+p_{4}\end{pmatrix} \\
  &= \begin{pmatrix}\left( x^{(1)}(1+q_{1}^{j})+q_{2}^{j} \right)(1+p_{1})+p_{2}\\ \left( x^{(2)}(1+q_{3}^{j})+q_{4}^{j} \right)(1+p_{3})+p_{4}\end{pmatrix},\\
  \overline{y}_{j}\left( \overline{\overline x}(p),q^{j} \right)&:=y_j\left(\overline{\overline x}(p)\right)(1+q_5^j)+q_6^j, \\
  \widetilde{\overline \phi}\left(\overline{\overline x}(p),p,a^{(\psi)} \right) &:= \widetilde{\phi }\left( \overline{\overline x}(p),a^{(\psi)}\right)+p_{5} \\
  &=\overline{y}\left( \overline{\overline x}(p),q\right)V_{l}\Sigma_{l}^{-1}a^{(\psi)}+p_{5}. \label{secGPOD:GPOD}
\end{align} 
Here $p$ does not contain a scaling parameter like $q^{j}$ in (\ref{secdatabase:ytransform}), because it would be linear dependent on $a^{(\psi)}$. The linear least squares problem (\ref{secGPOD:minsum}) is then augmented by the transformation parameters $p\in\R^{5}$ and again a penalty term is included:
\begin{equation}
  \min_{a^{(\psi)}\in\R^{l},p\in\R^{5}}\frac{1}{2}\sum_{i=1}^{n} \biggl( \phi(x_{i})- \widetilde{\overline \phi}\left(\overline{\overline x}_{i}(p),p,a^{(\psi)} \right)  \biggr)^{2}+\frac{\delta}{2}p^\top p. \label{secGPOD:mintransform2}
\end{equation}
This nonlinear least squares problem is solved by a Gau{\ss}-Newton algorithm for $a$ and $p$ simultaneously. Unlike (\ref{secdatabase:optproblem}), this can be accomplished without computational issues. Intuitively, an initial value for the algorithm is $p=0$ (no transformation) and its corresponding linear least squares solution $a^{(\psi)}=(\Psi^\top\Psi)^{-1}\Psi^\top\varphi$ (\ref{secGPOD:linsystem}). 

We call the approximation $\widetilde{\overline \phi}\left(\overline{\overline x}(p),p,a^{(\psi)} \right)$ to the unknown function $\phi(x)\in\LL$ based on the data pairs of evaluations $\left\{ (x_{i};\phi(x_{i})) \right\}_{i=1}^{n}$ and the function database $\left\{y_{1},\dots,y_{m}\right\}$ a \textit{generic surrogate model} (GSM). It is a least squares approximation to the data pairs and not an interpolation like the Kriging surrogate model. We consider the information contained in the data points as extremely valuable, especially since the evaluations of $\phi(x_{i})$ are assumed computationally very expensive. So a surrogate model should be as accurate as possible particularly in the proximity of any $x_{i}$, which can be realized by interpolation rather than approximation. Therefore, the next section will introduce a Kriging type data fusion method to generate an interpolation model which uses the generic surrogate model as a global trend. 

\section{Hierarchical Kriging}\label{sec:HK}

\textit{Variable fidelity modeling} (VFM) comprises methods for improving the approximation quality of interpolating only few computationally expensive evaluations (high fidelity) when having access to secondary data. This secondary data may consist of cheaper computations of a less accurate model (low fidelity) or a second variable which is assumed to be correlated with the primary variable. E.g.\ in CFD, computations with a Navier-Stokes code are regarded as high fidelity data, while computations of the same problem with an Euler code serve as low fidelity data. In this paper, we use the generic surrogate model $\widetilde{\overline \phi}$ (\ref{secGPOD:GPOD}) as the low fidelity model to improve the interpolation quality of the (high fidelity) evaluations $\left\{ (x_{i};\phi(x_{i})) \right\}_{i=1}^{n}$. So far, two major VFM frameworks could be distinguished. \textit{Cokriging}, originally developed in geostatistics, establishes a relation between primary and auxiliary variable by cross correlation \cite{Wackernagel03}. Other VFM methods use an (additive, multiplicative or hybrid) \textit{bridge function}, which corrects the discrepancy between a lo-fi and a hi-fi surrogate model \cite{HanGH2010}. Recently, a new robust VFM method was introduced \cite{HanGoertz2011}, whose implementation and computational complexity does not exceed the common Kriging method. It is called \textit{hierarchical Kriging} and the ansatz is a straightforward extension of section \ref{sec:rsm}.

In the Kriging model (\ref{sec1:KrigingModel}) 
\begin{equation}
  \phi(x)=\beta \widetilde{\overline \phi}\left(\overline{\overline x}(p),p,a \right)+z(x), \quad x\in\Omega\subset\R^{d}, \label{sec1:VFMmodel}
\end{equation}
the regression term is replaced by the lo-fi model (in our case the generic surrogate model $\widetilde{\overline\phi}$) with $\beta\in\R$ and $z(x)$ capturing the lack of fit like in (\ref{sec1:process}). Analogously  to (\ref{sec1:Kriging}), the hierarchical Kriging predictor 
\begin{equation}
  \widehat \phi(x)=c(x)^{\top}\varphi 
\end{equation}
is a weighted sum of $\varphi=\left(\phi(x_{1}),\dots,\phi(x_{n})\right)^{\top}$. Again, the weights $c_{i}(x)$ are determined by solving the linear equation
\begin{equation}
  \begin{bmatrix}R & \Phi \\\Phi\T & 0 \end{bmatrix} \begin{pmatrix} c(x) \\ \mu(x) \end{pmatrix} = \begin{pmatrix} r(x) \\ \widetilde{\overline \phi}(\overline{\overline x}(p),p,a )\end{pmatrix} 
\end{equation}
which minimizes the mean squared error subject to the unbiasedness constraint (\ref{sec1:BLUE}). $R$, $r(x)$, $c(x)$ and $\mu(x)$ are the same as in (\ref{sec1:krigingsystem}), while $\Phi=\left[\widetilde{\overline \phi}(\overline{\overline x}_{i}(p),p,a )\right]_{i}\in\R^{n\times 1}$ replaces the linear regression design matrix and $\widetilde{\overline \phi}(\overline{\overline x}(p),p,a )\in\R$ contains the evaluation of the generic surrogate model in $x$.
In the same manner we derive a closed form of the hierarchical Kriging predictor
\begin{equation}
  \widehat \phi(x)=c(x)^{\top}\varphi = \begin{pmatrix} r(x) \\ \widetilde{\overline \phi}(\overline{\overline x}(p),p,a )\end{pmatrix}^{\top} \begin{bmatrix}R & \Phi \\\Phi\T & 0 \end{bmatrix}^{-1}\begin{pmatrix} \varphi \\ 0\end{pmatrix} \label{secHK:HK1}
\end{equation}
or equivalently with $\beta=\left(\Phi\T R^{-1}\Phi \right)^{-1}\Phi\T R^{-1}\varphi$ 
\begin{equation}
  \widehat \phi(x)=\beta \widetilde{\overline \phi}(\overline{\overline x}(p),p,a ) +r(x)\T R^{-1}\left( \varphi - \beta \Phi\right).\label{secHK:HK2}
\end{equation}
Formulation (\ref{secHK:HK1}) is more suitable for the implementation, because $\beta$ does not need to be computed explicitly and the solution of the linear system $ \left[\begin{smallmatrix}R & \Phi \\\Phi\T & 0 \end{smallmatrix}\right]^{-1}\left(\begin{smallmatrix} \varphi \\ 0\end{smallmatrix}\right)$ only has to be computed once since it is independent of $x$. Also one can witness that, inserting any $x_{i}$ $(i=1,\dots,n)$, the hierarchical Kriging predictor is an interpolator ($\widehat\phi(x_{i})=\phi(x_{i})$). Formulation (\ref{secHK:HK2}) demonstrates how the predictor works. Assuming that the generic surrogate model fit $\widetilde{\overline\phi}$ approximates the data $\phi(x_{i})$, $\beta$ will be close to $1$. The second summand is a weighted sum of correlation functions which ``pulls'' the response towards the exact evaluations $\phi(x_{i})$.

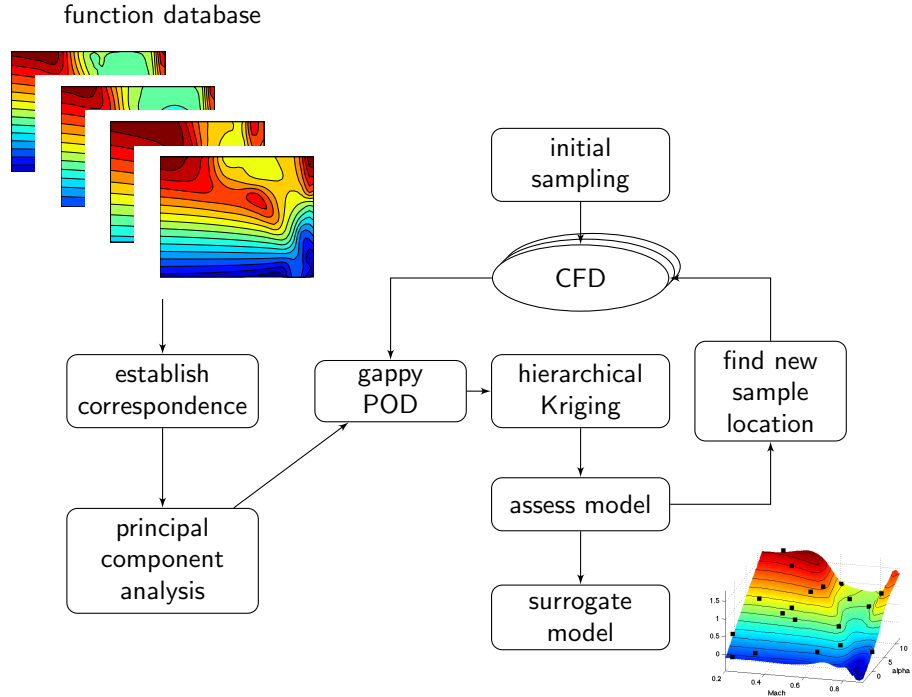
\begin{figure}[h!]
\centering{ 
\begin{small}
 
  \tikzstyle{block} = [rectangle, draw, text width=6em, text centered, rounded corners, minimum height=2em,font=\sffamily]
  \tikzstyle{block_new} = [rectangle, draw, text width=6.5em, text centered, rounded corners, minimum height=2em,font=\sffamily]
  \tikzstyle{line} = [draw, -latex']
  \tikzstyle{cloud} = [draw, ellipse,text centered,text width=4em, minimum height=2.5em,font=\sffamily]
    
  \begin{tikzpicture}[node distance = 1.5cm,auto]
    
    \node [block] (Initial DoE) {initial sampling};
    \node [cloud, double copy shadow, fill=white, below of=Initial DoE] (CFD) {CFD};
    \node [block, below of=CFD] (Kriging) {hierarchical\\Kriging};
    \node [block, below of=Kriging] (Assess model) {assess model};
    \node [block, below of=Assess model] (Surrogate model) {surrogate model};
    \node [block, right of=Kriging, node distance=2.5cm,text width=5em] (new sample) {find new sample location};
   
    \node [right of=Surrogate model,node distance=3.0cm] {\includegraphics[width=0.2\textwidth]{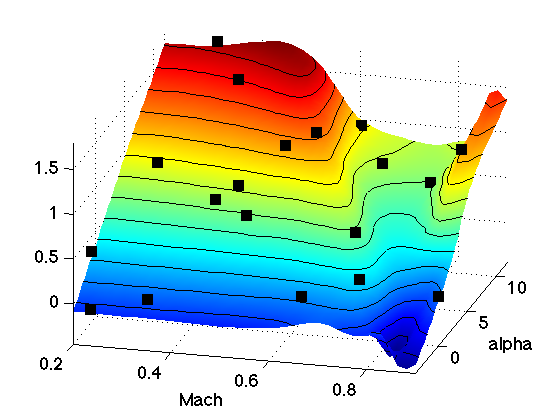}};

    \node [block_new, left of=Kriging, node distance=2.5cm,text width=5em] (generic model) {gappy\\ POD};

    \node [block_new, left of=generic model, node distance=3cm] (Fitting) {establish\\correspondence};
    \node [block_new, below of=Fitting,node distance=2.25cm] (PCA) {principal component analysis};
    \node [above of=Fitting, node distance=3cm,text centered  ] (database) {\includegraphics[width=0.3\textwidth]{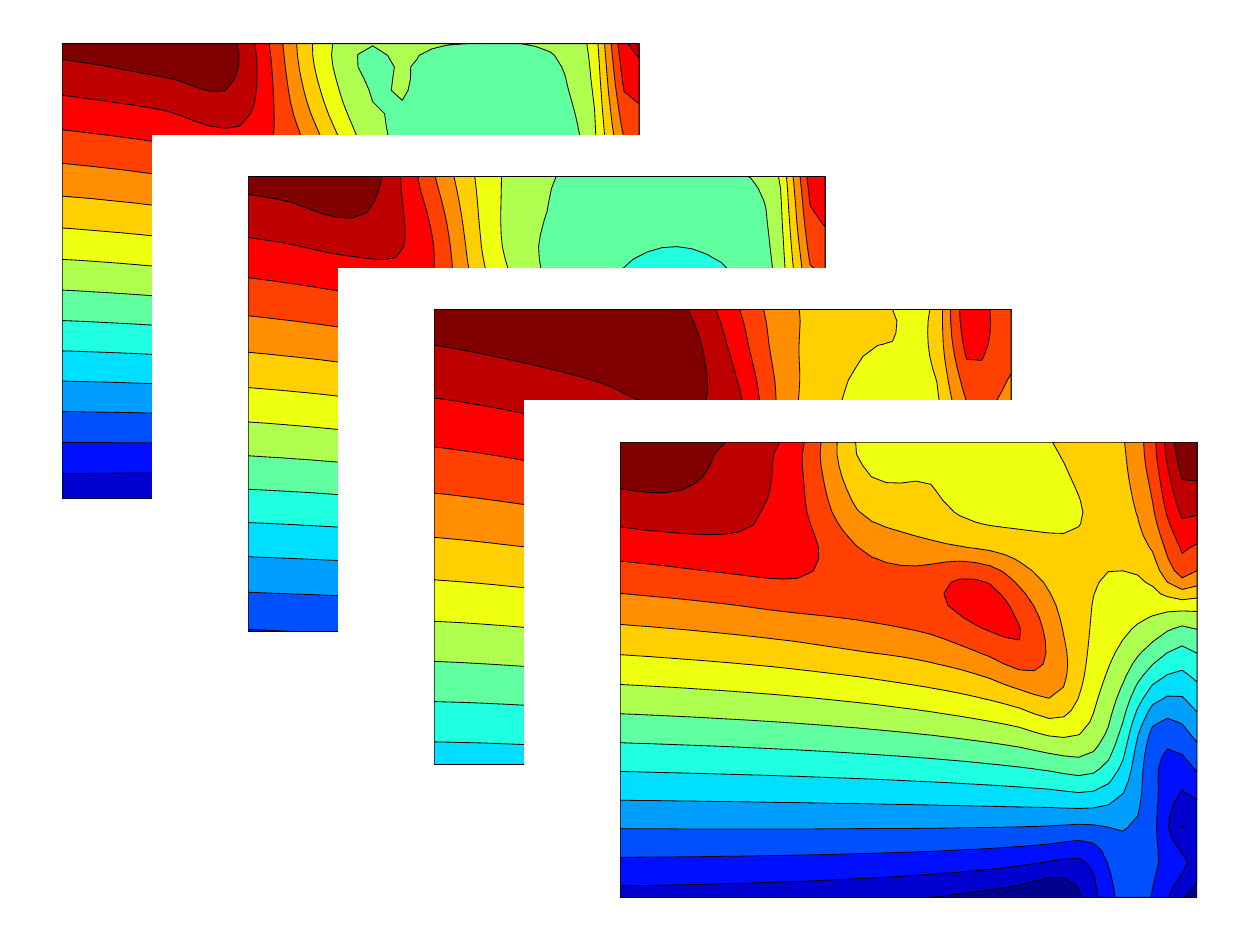}  };
    \node [above of=database, node distance=2.0cm,text centered,font=\sffamily ] (database_text) { function database};
  
    \path [line] (Initial DoE) -- (CFD);
    \path [line] (Kriging) -- (Assess model);
    \path [line] (Assess model) -- (Surrogate model);
    \path [line] (Assess model) -| (new sample);
    \path [line] (new sample) |- (CFD);

    \path [line] (Fitting) -- (PCA);
    \path [line] (database) -- (Fitting);

    \path [line] (PCA) -- (generic model);
    \path [line] (CFD) -| (generic model);
    \path [line] (generic model) -- (Kriging);

  \end{tikzpicture}
\end{small}
}
\caption{Generic surrogate modeling framework}
\label{fig:framework}
\end{figure}

Figure \ref{fig:framework} outlines the generic surrogate modeling framework. We distinguish between offline costs and online costs for model generation. The three steps on the left hand side, generation of accurate surrogate models for all database functions, solving the correspondence problem and computing the proper orthogonal decomposition, only have to be performed once so they are merely offline costs. When interpolating sample data for any new test case, evaluating the computer experiment for the samples, the gappy POD fit and building the hierarchical Kriging model are online costs. We assume that one single evaluation of the computer experiment $\phi(x_{i})$, e.g.\ a CFD solution, dominates the computational cost of model generation.

\section{Numerical Results}\label{sec:num}

\begin{figure}[h]
\centering{ 
  \includegraphics[width=0.325\textwidth]{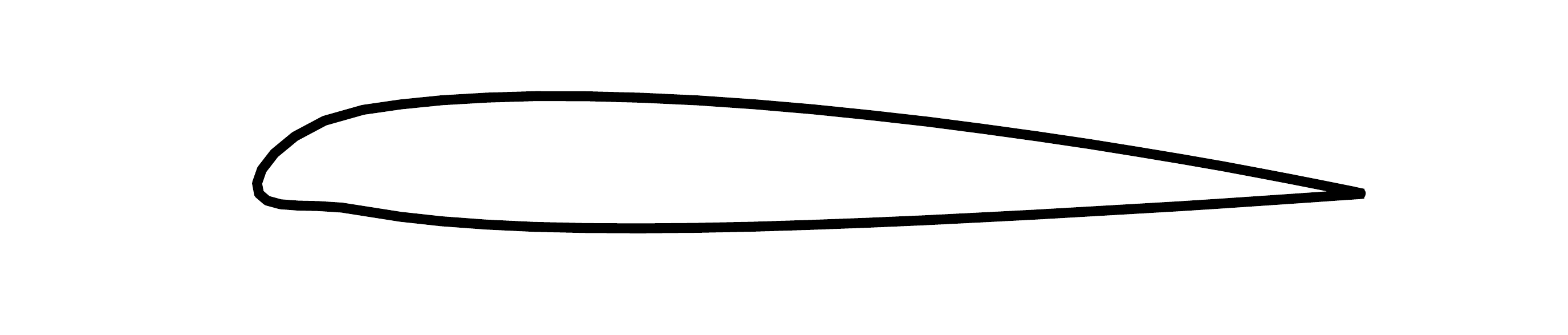}
  \includegraphics[width=0.325\textwidth]{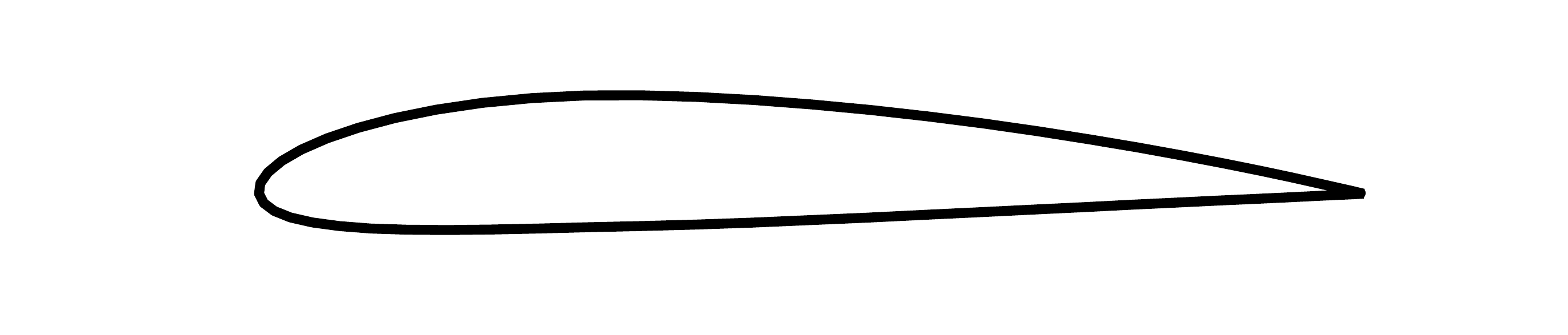}
  \includegraphics[width=0.325\textwidth]{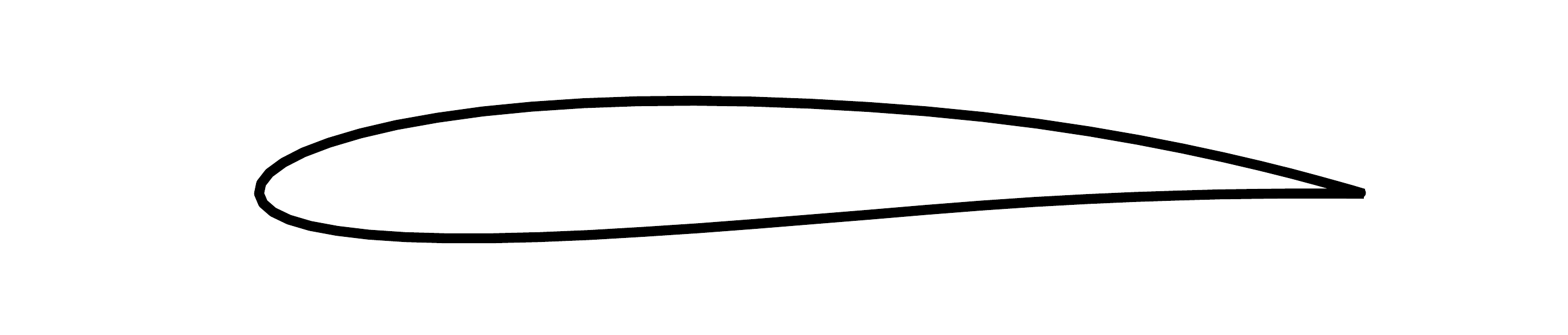}\\
  \includegraphics[width=0.325\textwidth]{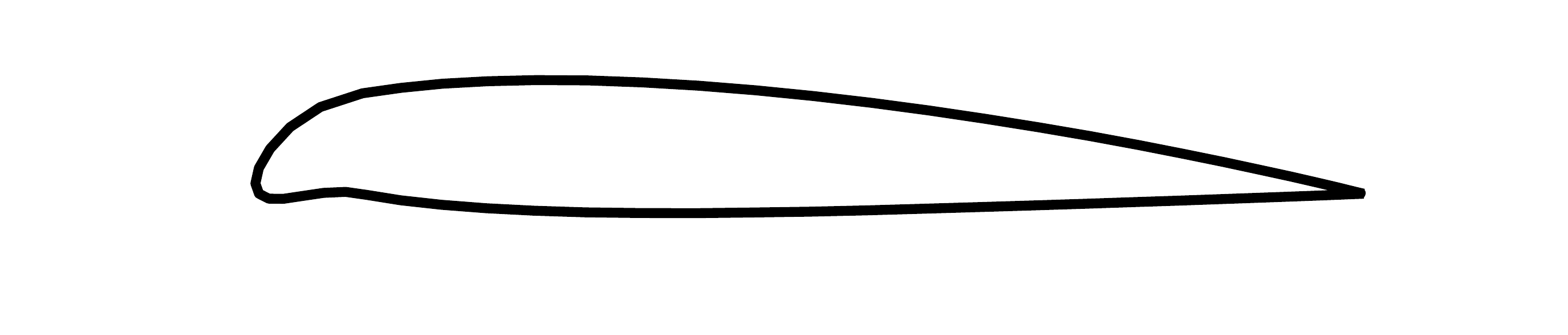}
  \includegraphics[width=0.325\textwidth]{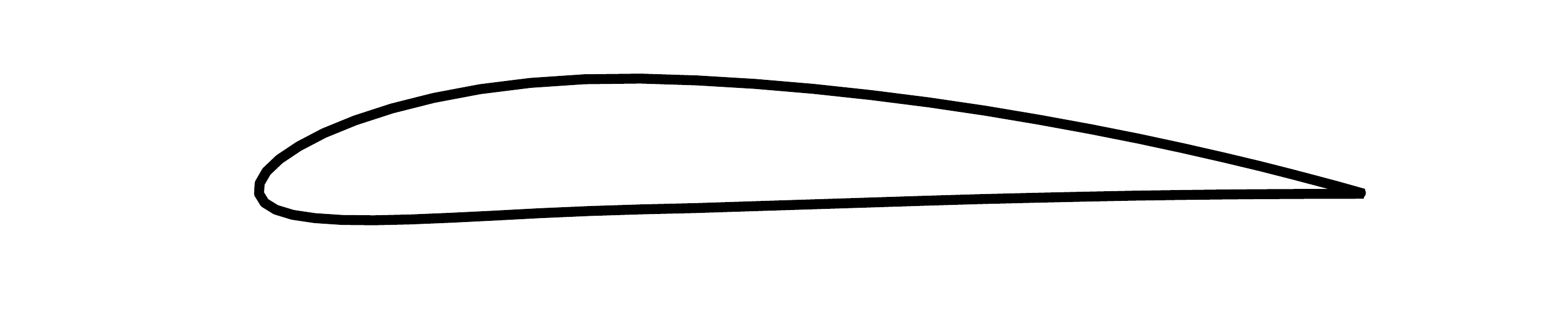}
  \includegraphics[width=0.325\textwidth]{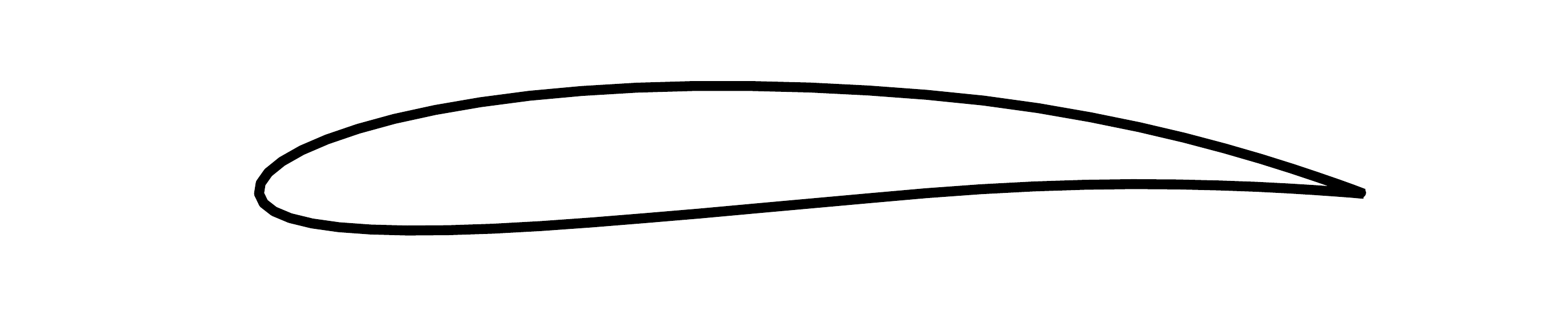}\\
  \includegraphics[width=0.325\textwidth]{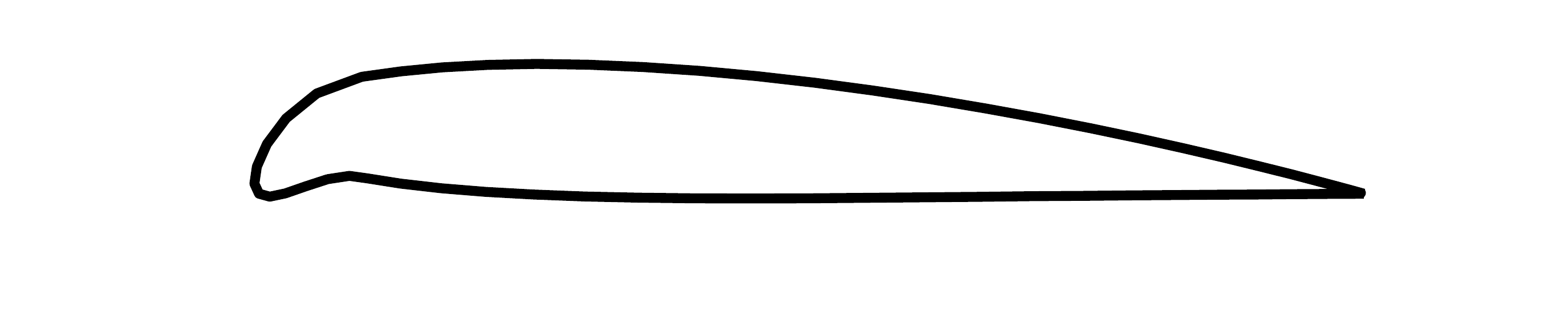}
  \includegraphics[width=0.325\textwidth]{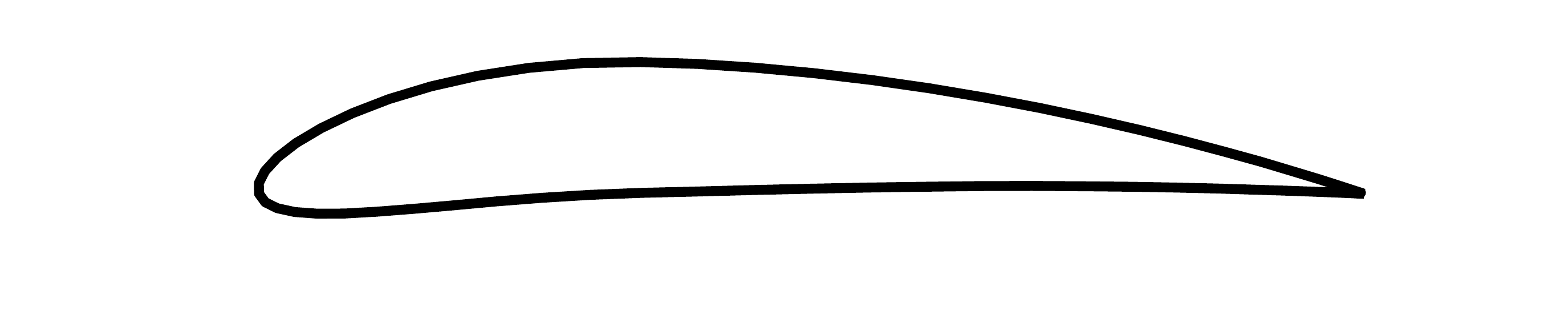}
  \includegraphics[width=0.325\textwidth]{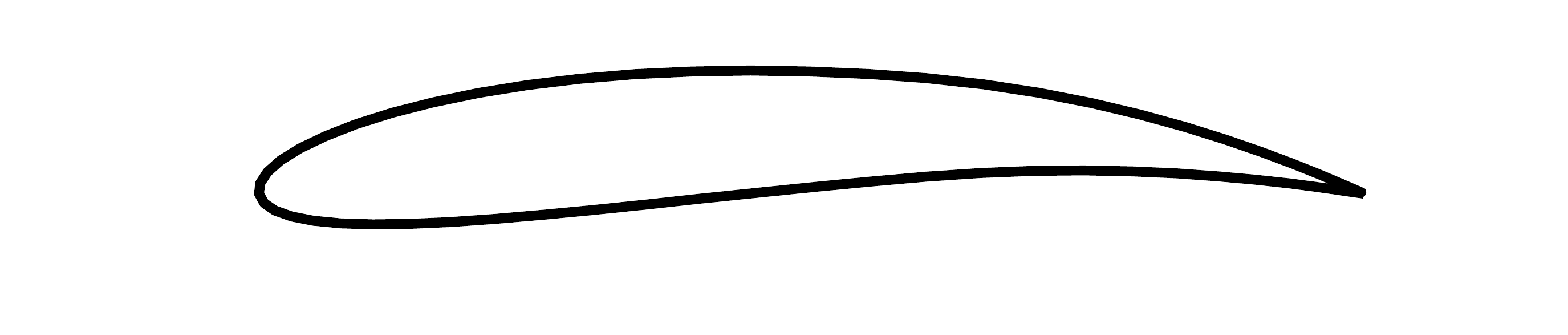}\\ 
  \includegraphics[width=0.325\textwidth]{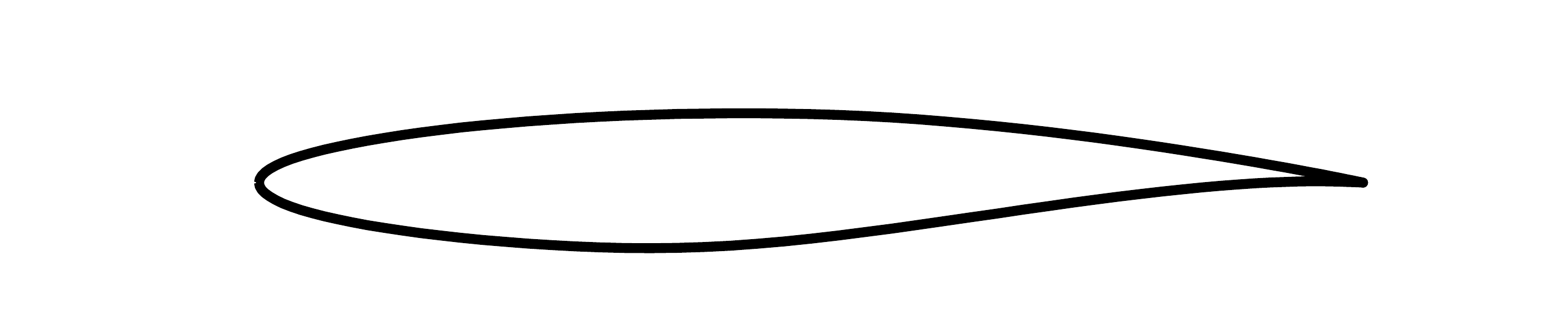}
}
\caption{Airfoil database, 9 examples from the NACA 4-digit series and the RAE 2882.}
\label{fig:airfoils}
\end{figure}

\begin{figure}[h]
\centering{ 
  \includegraphics[width=0.325\textwidth]{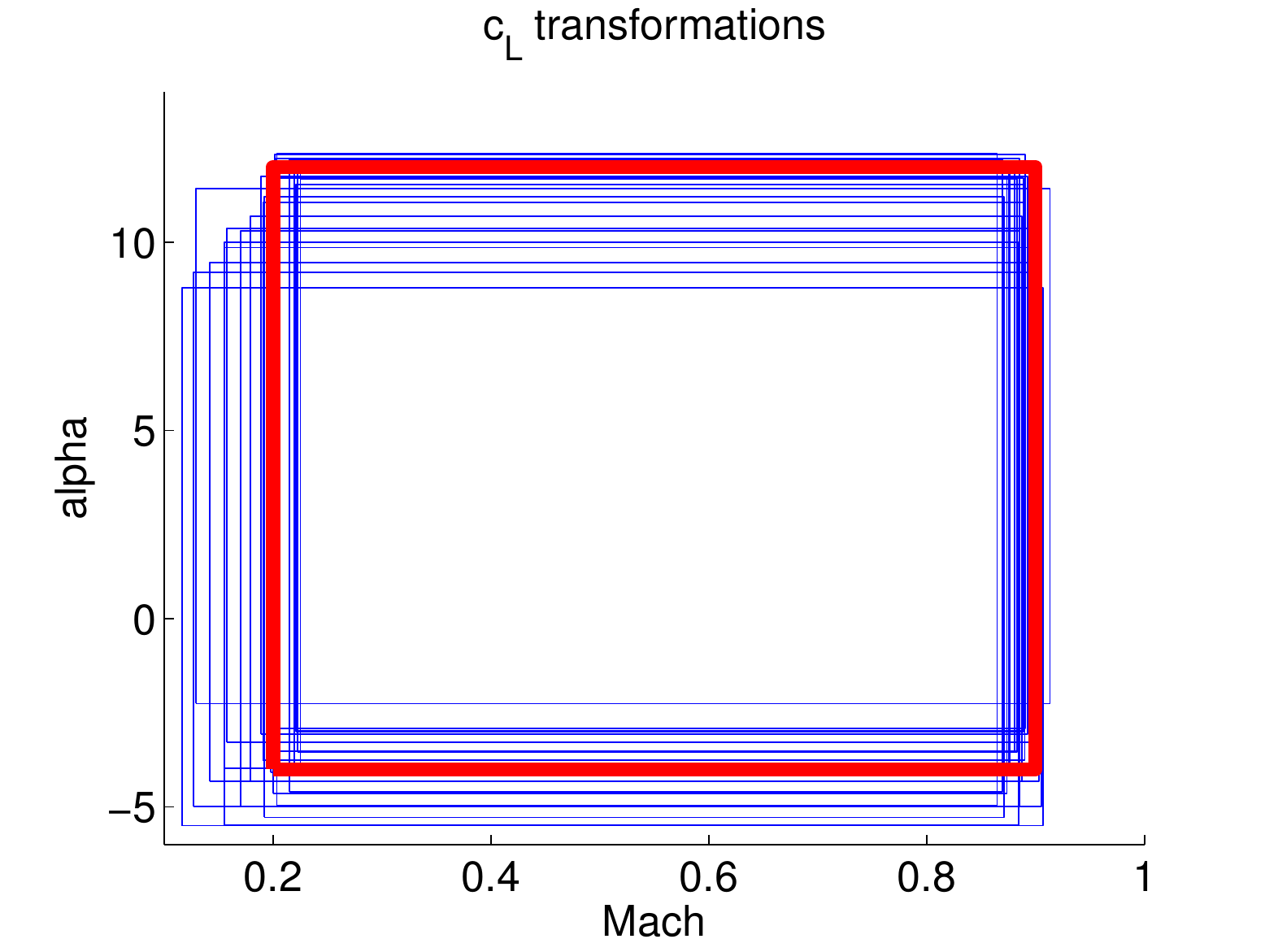}
  \includegraphics[width=0.325\textwidth]{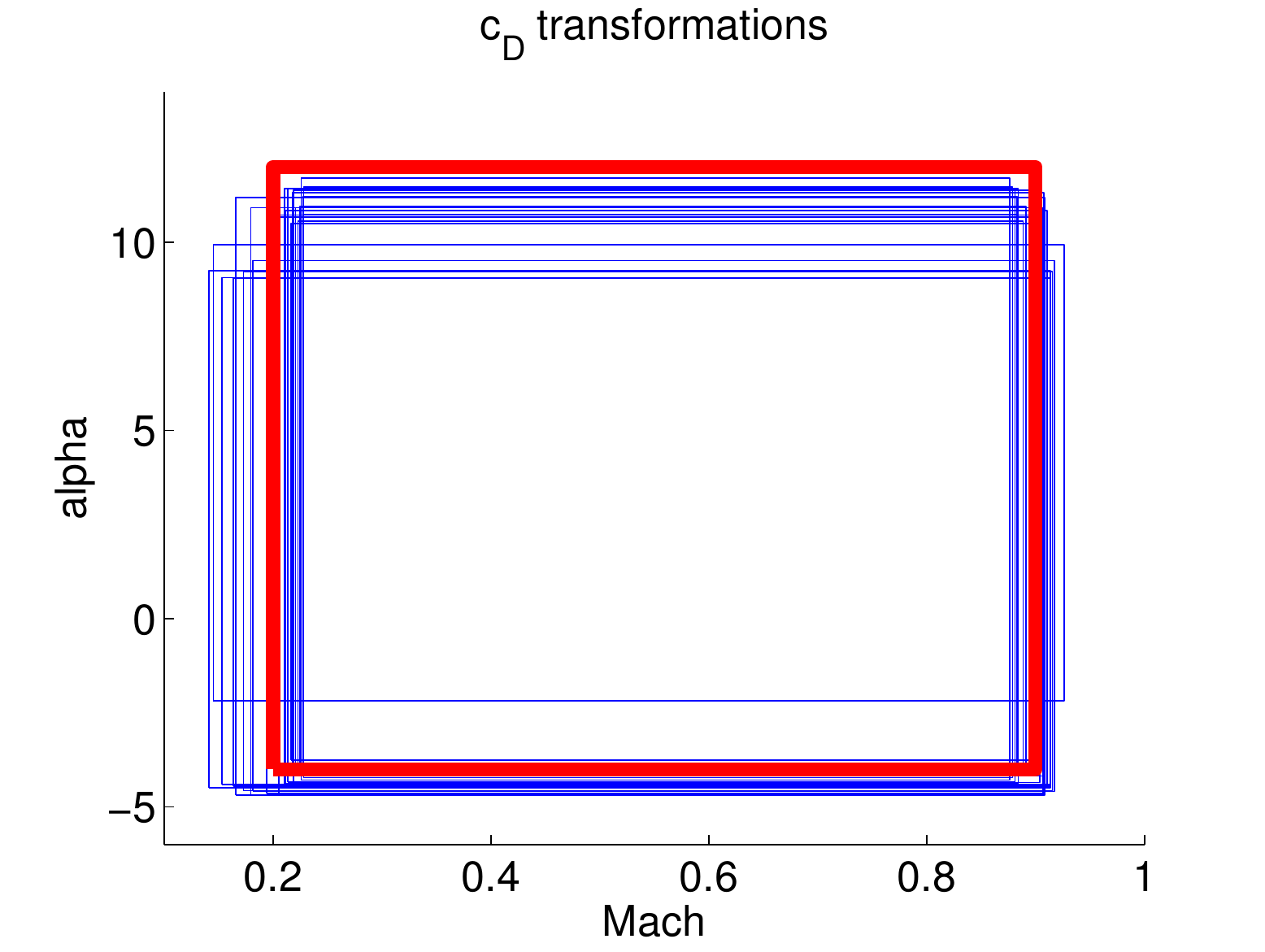}
  \includegraphics[width=0.325\textwidth]{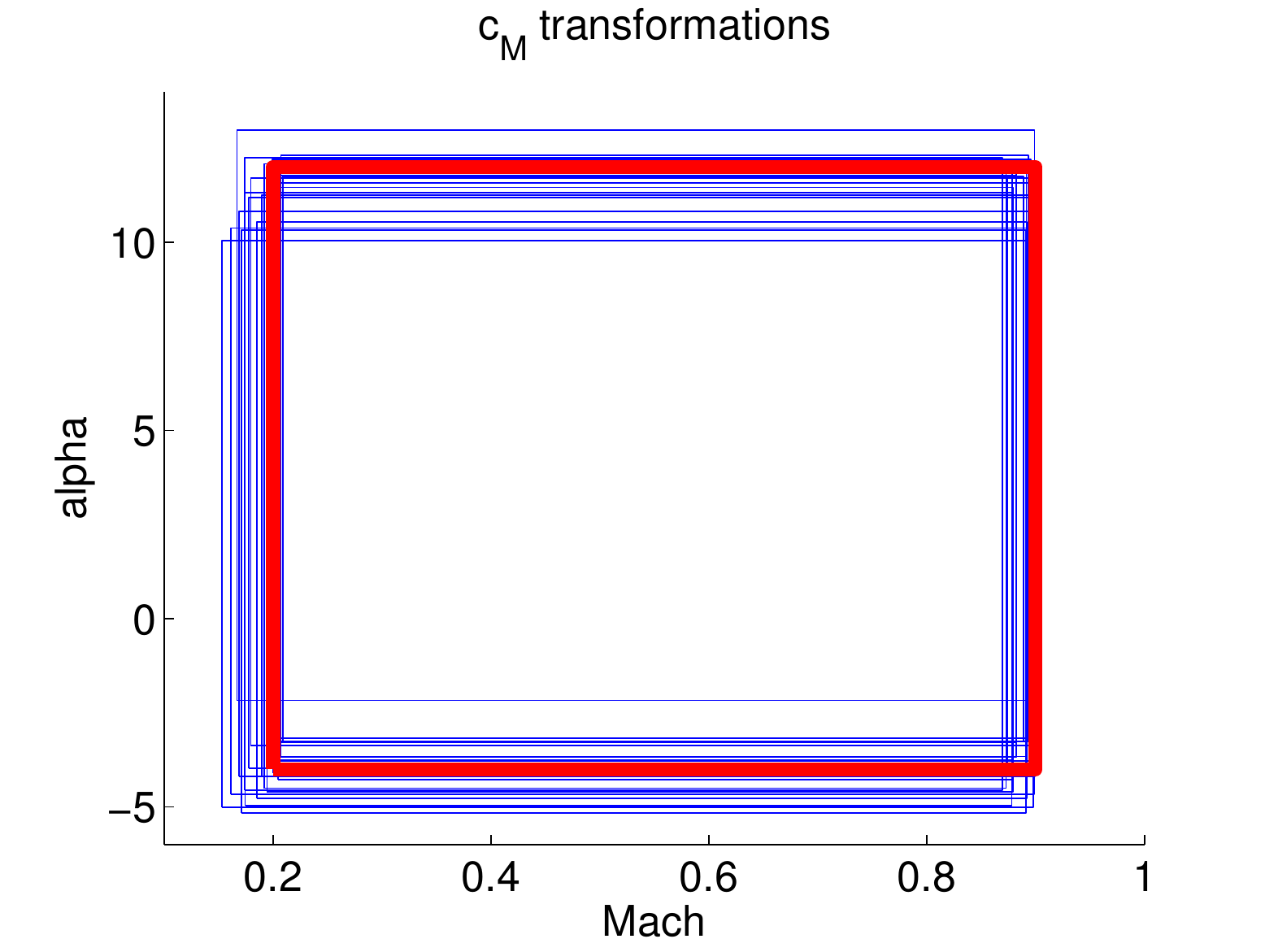}}
\caption{Red: reference domain $\Omega$, blue: preimages of $\Omega$ for the transformations (\ref{secdatabase:xtransform}).}
\label{fig:transformations}
\end{figure}

We choose the following test case for the validation of the methods of this paper. For two-dimensional airfoil geometries, we approximate the aerodynamic coefficients lift ($c_{l}$), drag ($c_{d}$) and pitching moment ($c_{m}$), depending on the Mach number and the angle of attack $\alpha$. The data is obtained by RANS computations with the DLR TAU-code \cite{taucode2006}. We consider a database of 24 airfoils: 23 airfoils from the NACA 4-digit series (first digit $\in \left[3.0,6.0\right]$, second digit $\in \left[1.0,6.0\right]$ and the last two digits which determine the thickness are fixed at $12$) and the RAE 2822, see figure \ref{fig:airfoils}. For each airfoil, we compute an accurate surrogate model for each problem class ($c_{l}(\textit{Ma},\alpha)$, $c_{d}(\textit{Ma},\alpha)$, $c_{m}(\textit{Ma},\alpha)$) in the reference domain $\Omega:=\left[0.2,0.9\right]\times\left[-4^{\circ},+12^{\circ} \right]$. Note that if alignment of the database is considered (\ref{secdatabase:optproblem}), the surrogates should be accurate in a domain larger than the reference domain (e.g.\ $\left[0.1,1.0\right]\times\left[-6^{\circ},+14^{\circ} \right]$) due to possible translations (\ref{secdatabase:xtransform}) in the alignment process, see figure \ref{fig:transformations}. For each airfoil, $400$ CFD solutions are computed on a $20\times 20$ tensorgrid which covers the reference domain. This sums up to 9600 CFD simulations for generating the database. Depending on the input parameters $(\textit{Ma},\alpha)$, one flow solution takes from 30 minutes to over 3 hours of CPU time. All computations can be performed independently from each other in parallel and using two AMD Opteron architectures with 48 2.3GHz cores each, the generation of the databases took approximately two weeks.

\begin{figure}[h!]
\centering{ 
  \includegraphics[width=0.49\textwidth]{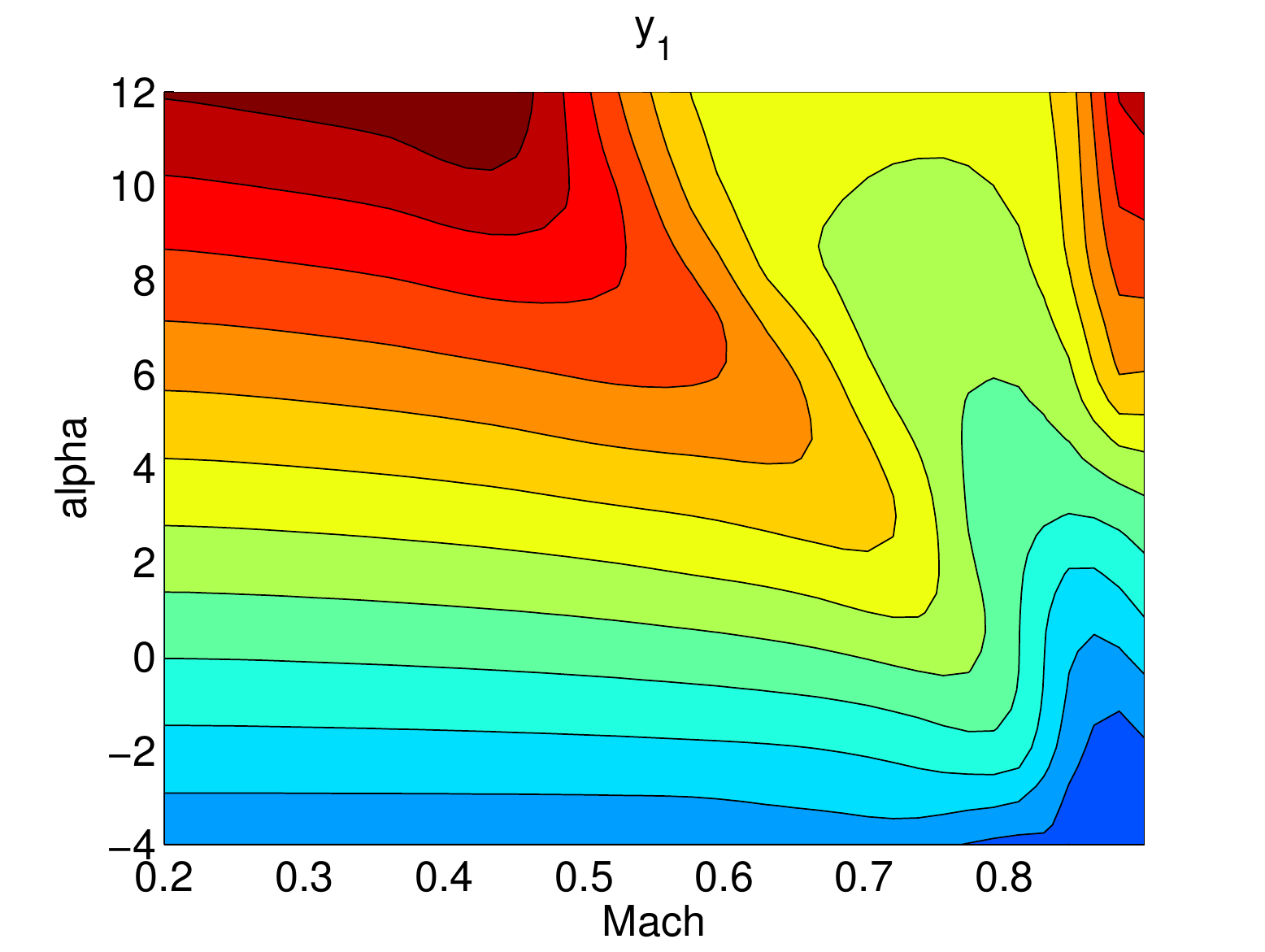}
  \includegraphics[width=0.49\textwidth]{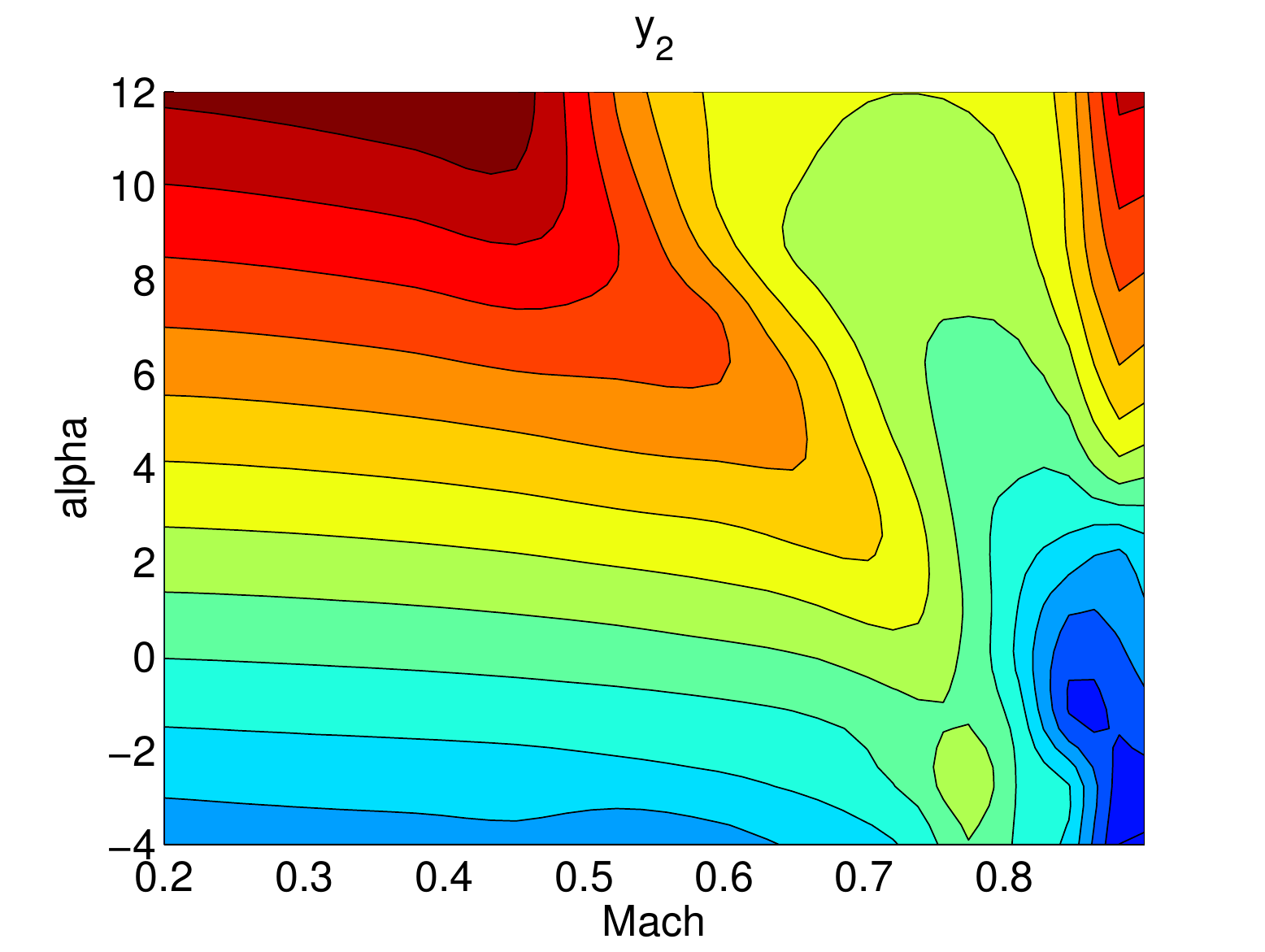}\\
  \includegraphics[width=0.49\textwidth]{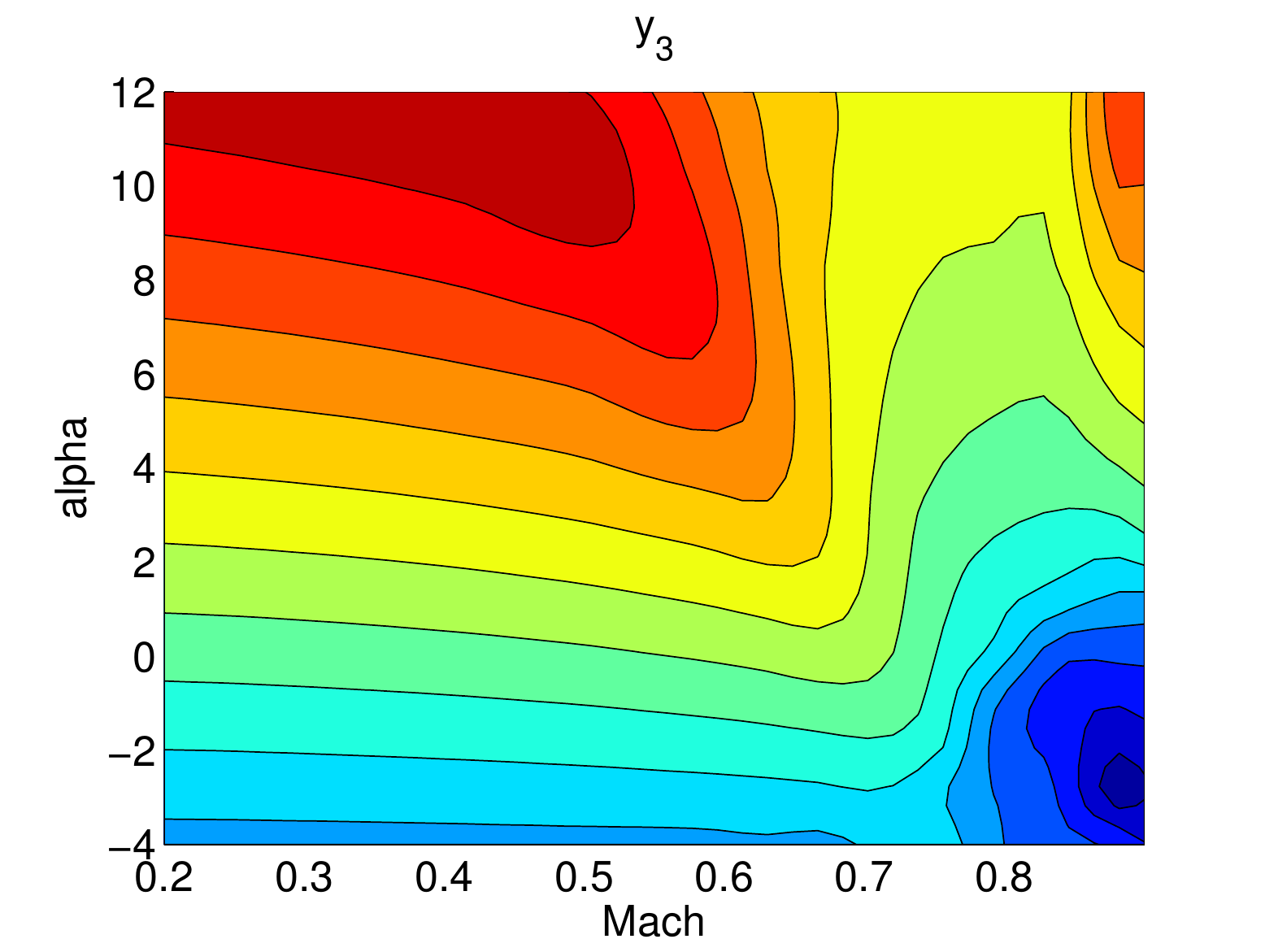}
  \includegraphics[width=0.49\textwidth]{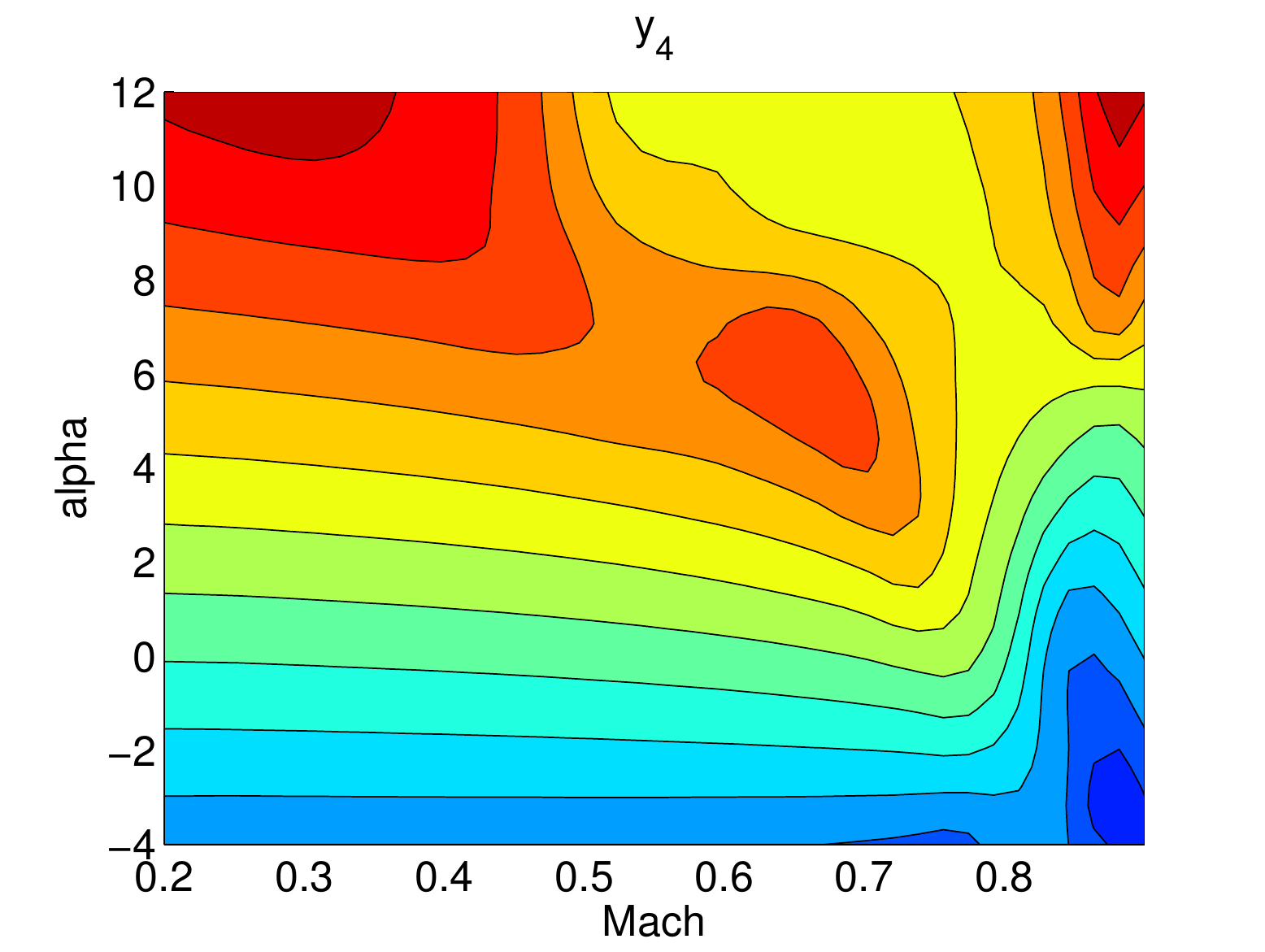}
}
\caption{$c_{l}$ database functions}
\label{fig:cldatabase}
\end{figure}

For each of the three test cases, after solving the correspondence problem (\ref{secdatabase:optproblem}) the database contains the aligned surrogate models $y_{1}(x),\dots,y_{24}(x)$. Figure \ref{fig:cldatabase} shows four examples of lift responses. Mutual characteristics like the linear behavior in the left part, the local maximum in the upper right corner or the curvature which reaches from the lower right to the upper middle are aligned as good as possible by the admissible transformation. 

\begin{figure}[h]
\centering{ 
  \includegraphics[width=0.48\textwidth]{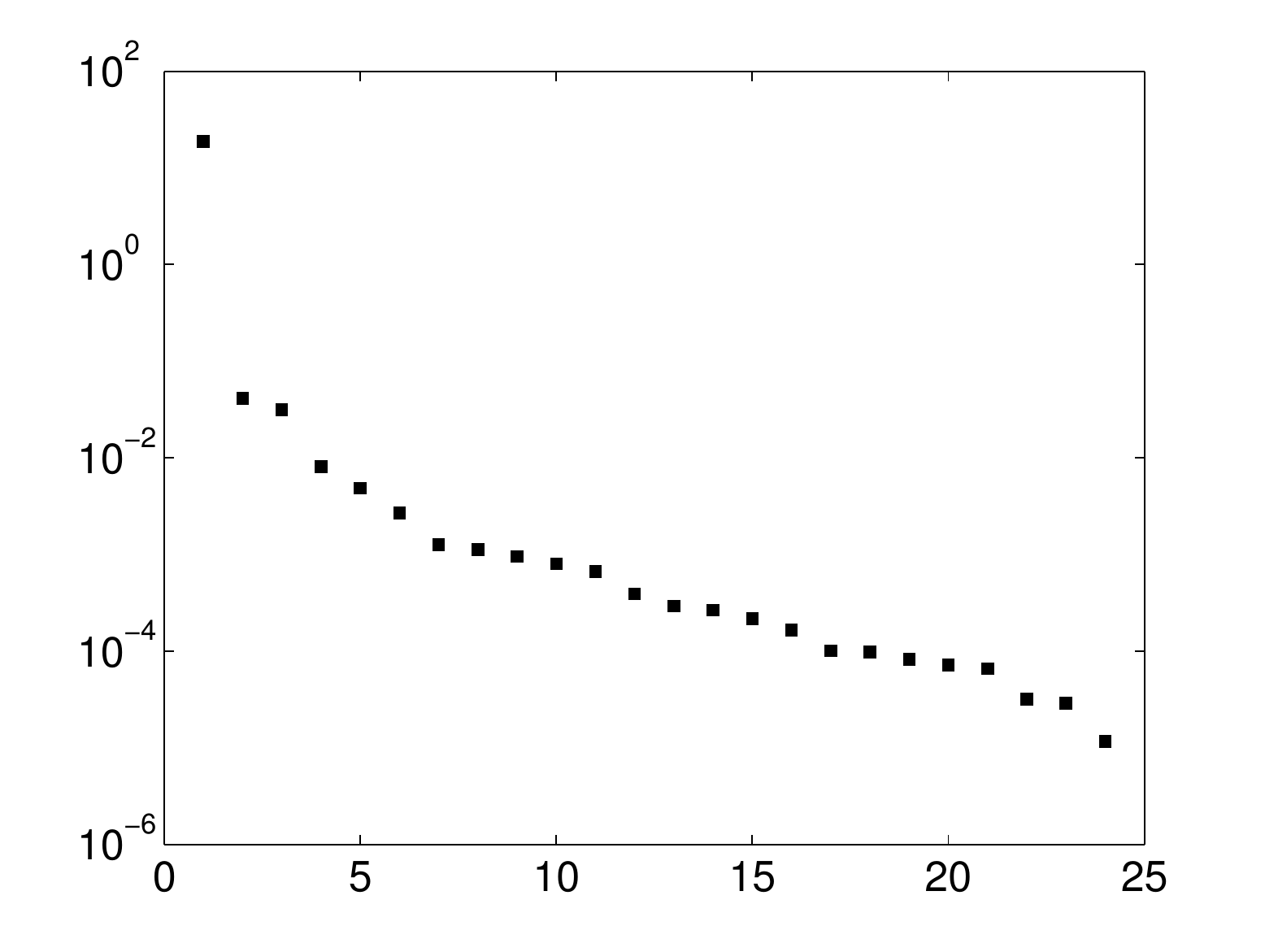} 
}
\caption{Distribution of the eigenvalues $\lambda_{i}$ for the $c_{l}$ POD basis}
\label{fig:eigenvalues}
\end{figure}

\begin{figure}[h!]
\centering{ 
  \includegraphics[width=0.49\textwidth]{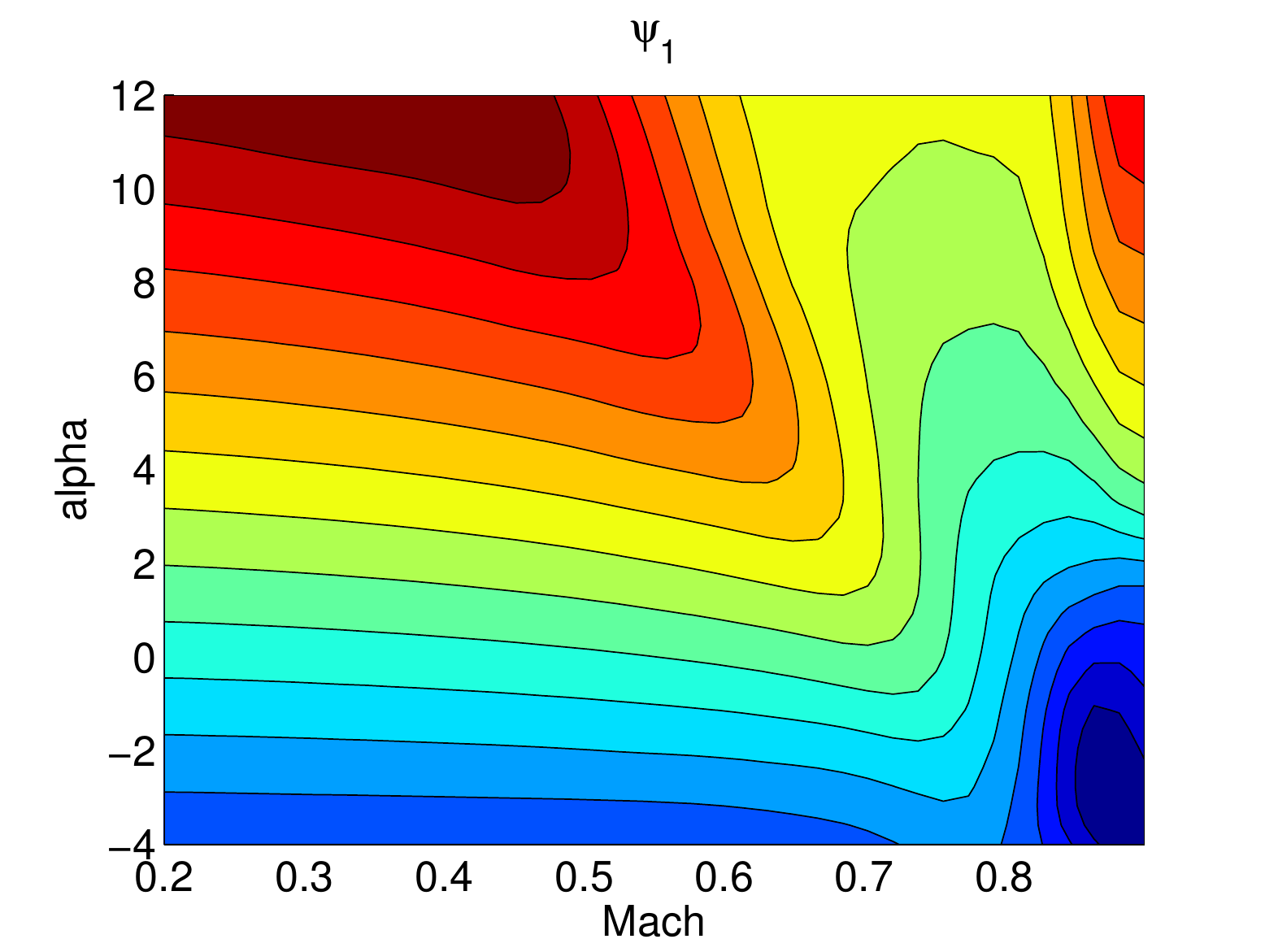}
  \includegraphics[width=0.49\textwidth]{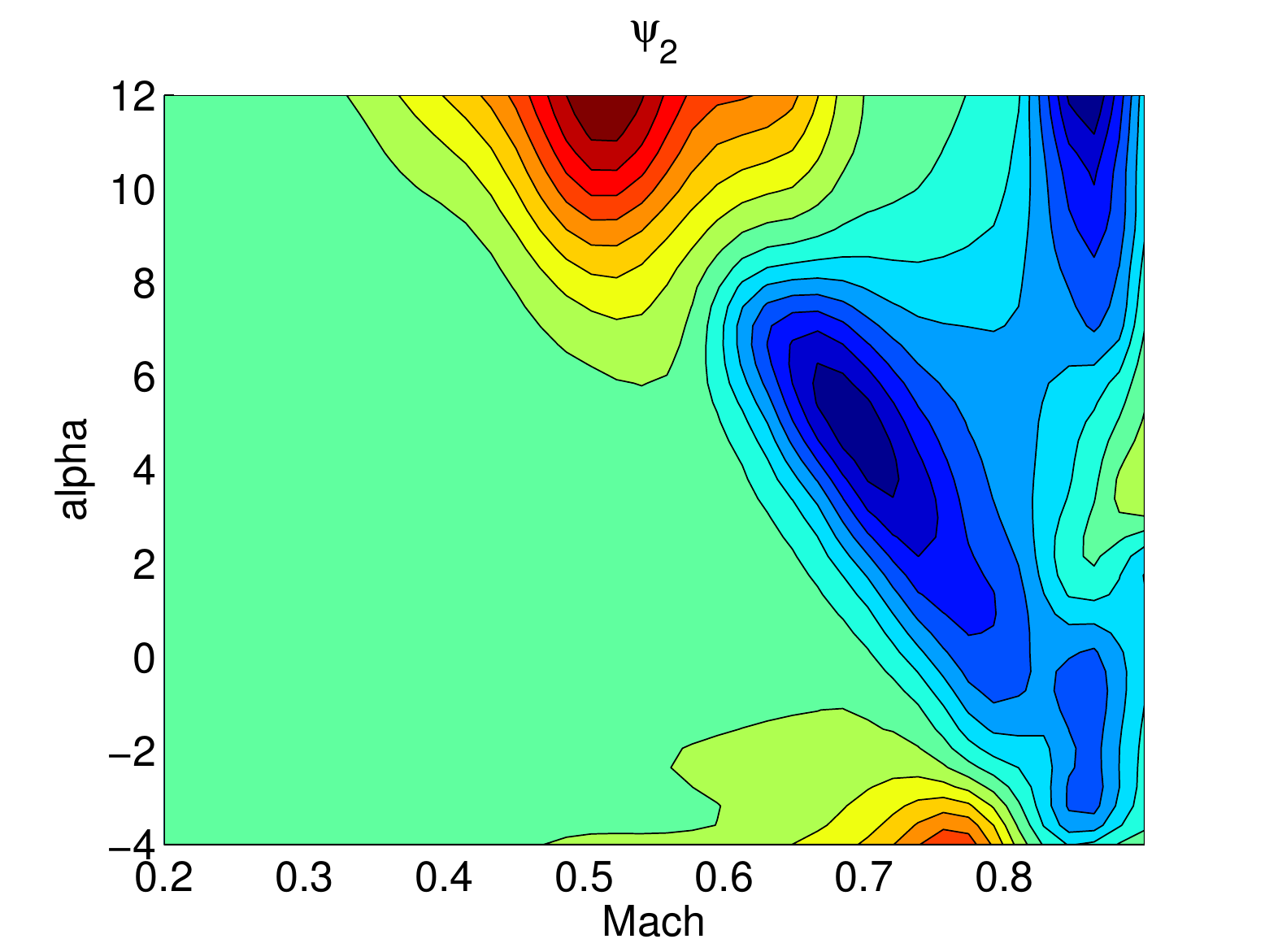}\\
  \includegraphics[width=0.49\textwidth]{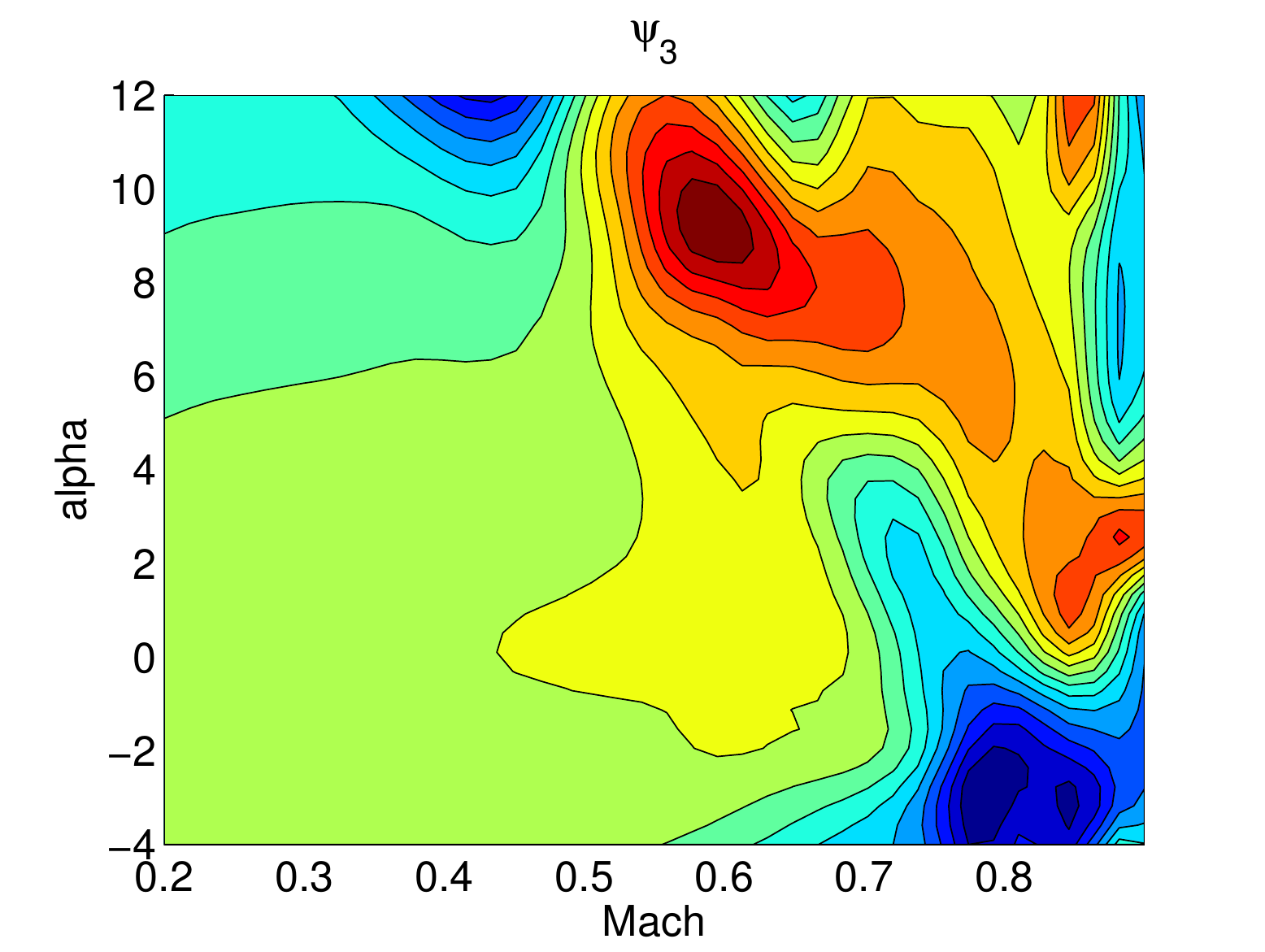}
  \includegraphics[width=0.49\textwidth]{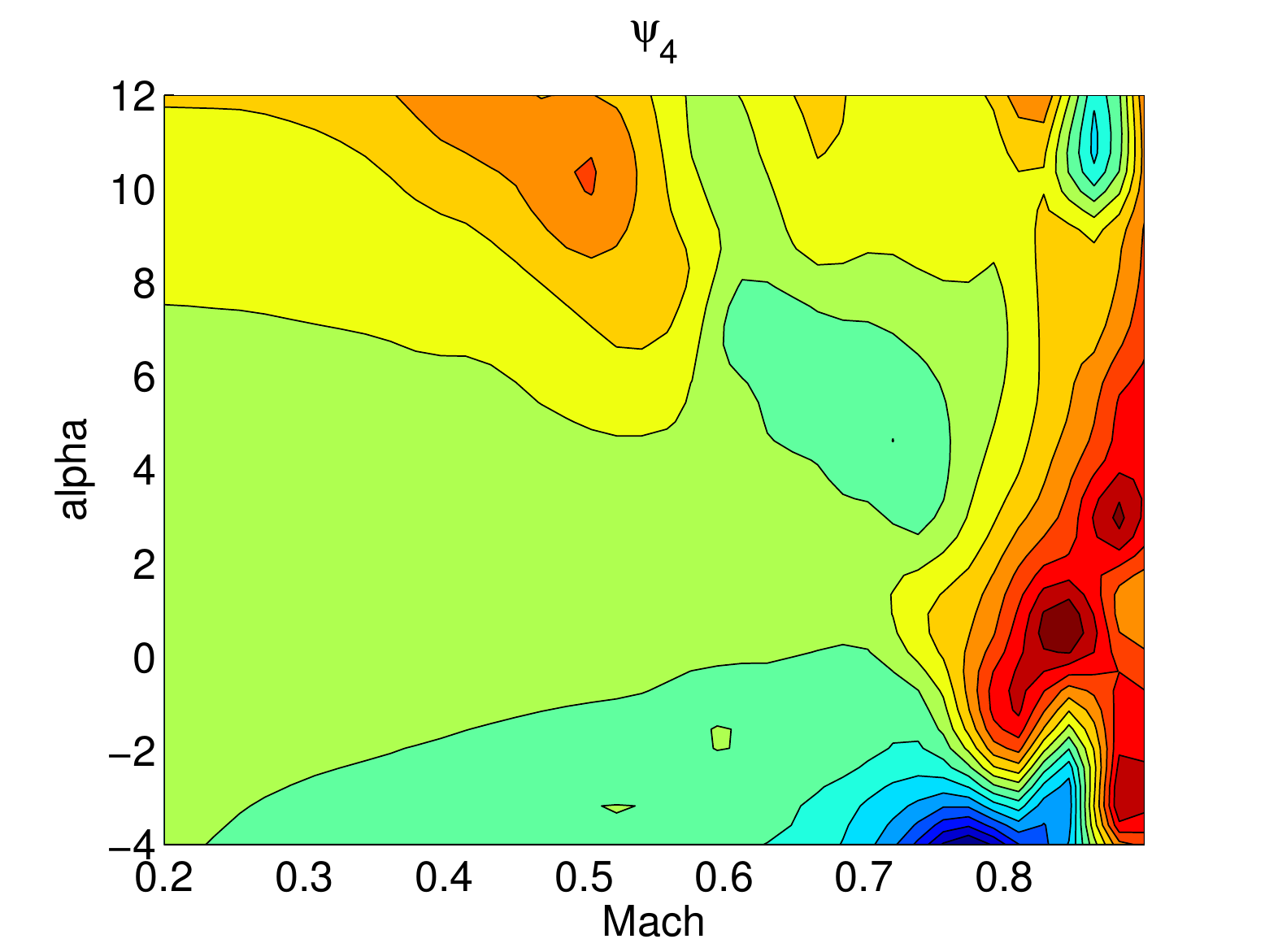}
}
\caption{$c_{l}$ POD basis}
\label{fig:clpodbasis}
\end{figure}

Subsequently, the proper orthogonal decomposition of $\mathcal{Y}=\operatorname{span}\left\{y_{1},\dots,y_{24}\right\}$ is performed. The choice of the POD-basis' rank $l$, i.e.\ the number of POD-basis elements, is carried out automatically. Using the approximation error formula (\ref{secPOD:errorformula}), $l$ is chosen such that 
\begin{equation}
  \frac{\sum_{j=1}^{l}\lambda_{j}}{\sum_{j=1}^{24}\lambda_{j}}\geq 0.999. \label{secnum:error}
\end{equation}
In this way, we guarantee that the POD-basis $\left\{\psi_{1},\dots,\psi_{l}\right\}$ contains at least $99.9\%$ of $\mathcal Y$'s information. Table \ref{tab:rank} shows the number of required POD-basis elements for the three test cases. Using an aligned database reduces the number by one in each case, meaning that by aligning the characteristic features of the database functions the total variation is reduced. Figure \ref{fig:eigenvalues} illustrates the rapid decay of the eigenvalues for the $c_{l}$ test case, therefore already four POD-basis elements are sufficient for (\ref{secnum:error}). The properties of $\psi_{1},\dots,\psi_{4}$ are depicted in figure \ref{fig:clpodbasis}: $\psi_{1}$ is close to an average $c_{l}$ response, $\psi_{2}$ and $\psi_{3}$ both essentially control the curvature from the lower right to the upper middle as well as the position of the local maximum in the upper right and $\psi_{4}$ mostly influences the behavior in the lower right.
\begin{table}[h]
\footnotesize\centering{
  \begin{tabular}{|l|ccc|}\firsthline
     & $c_{l}$ & $c_{d}$ & $c_{m}$ \\ \hline
    POD (aligned) & 4 & 4 & 5 \\ 
    POD (no alignment)  & 5 & 5 & 6\\ \lasthline
  \end{tabular}
}
  \caption{Rank of the POD-basis.}
  \label{tab:rank}
\end{table}

\begin{figure}[h!]
\centering{ 
  \includegraphics[width=0.49\textwidth]{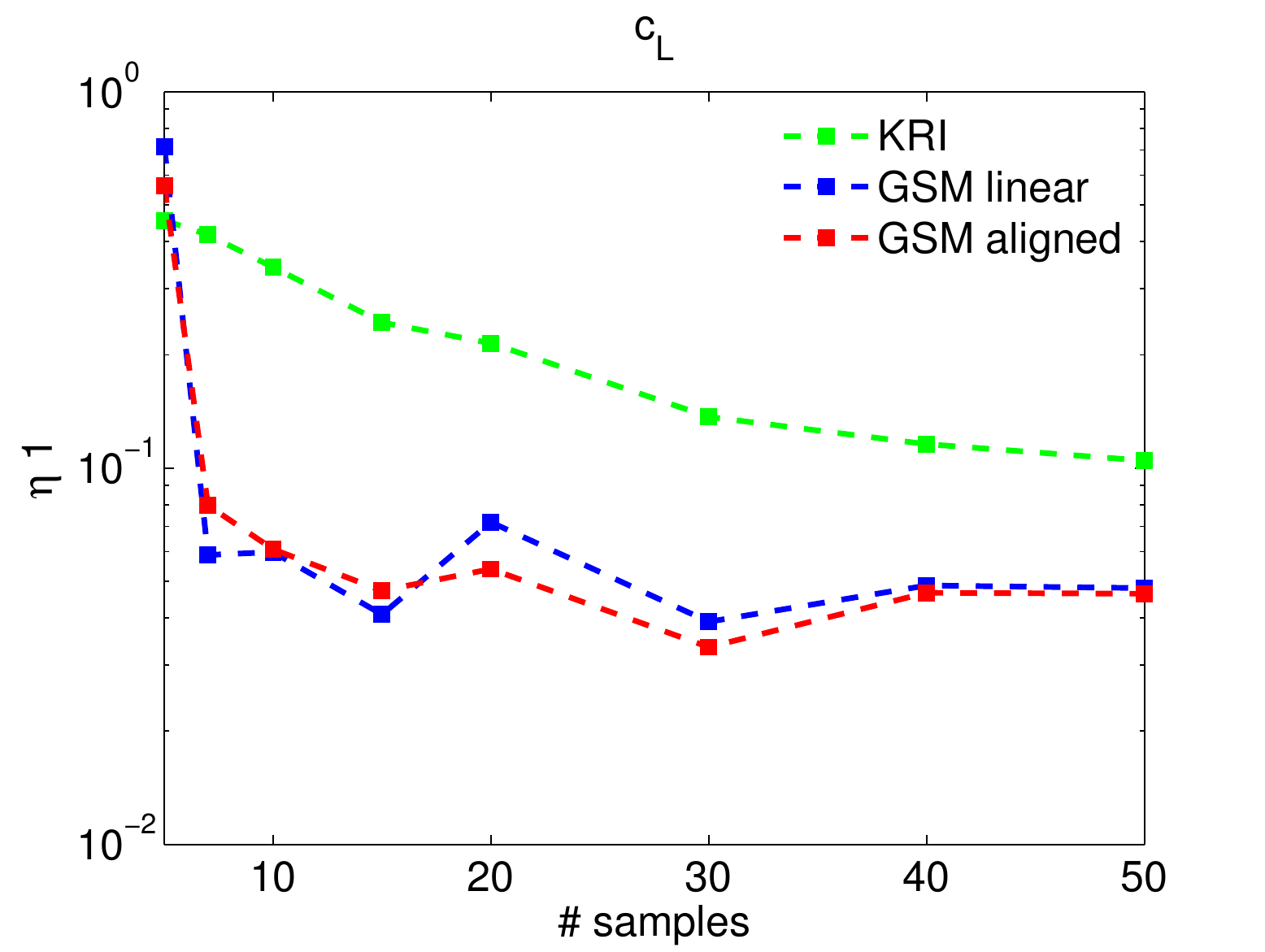} 
  \includegraphics[width=0.49\textwidth]{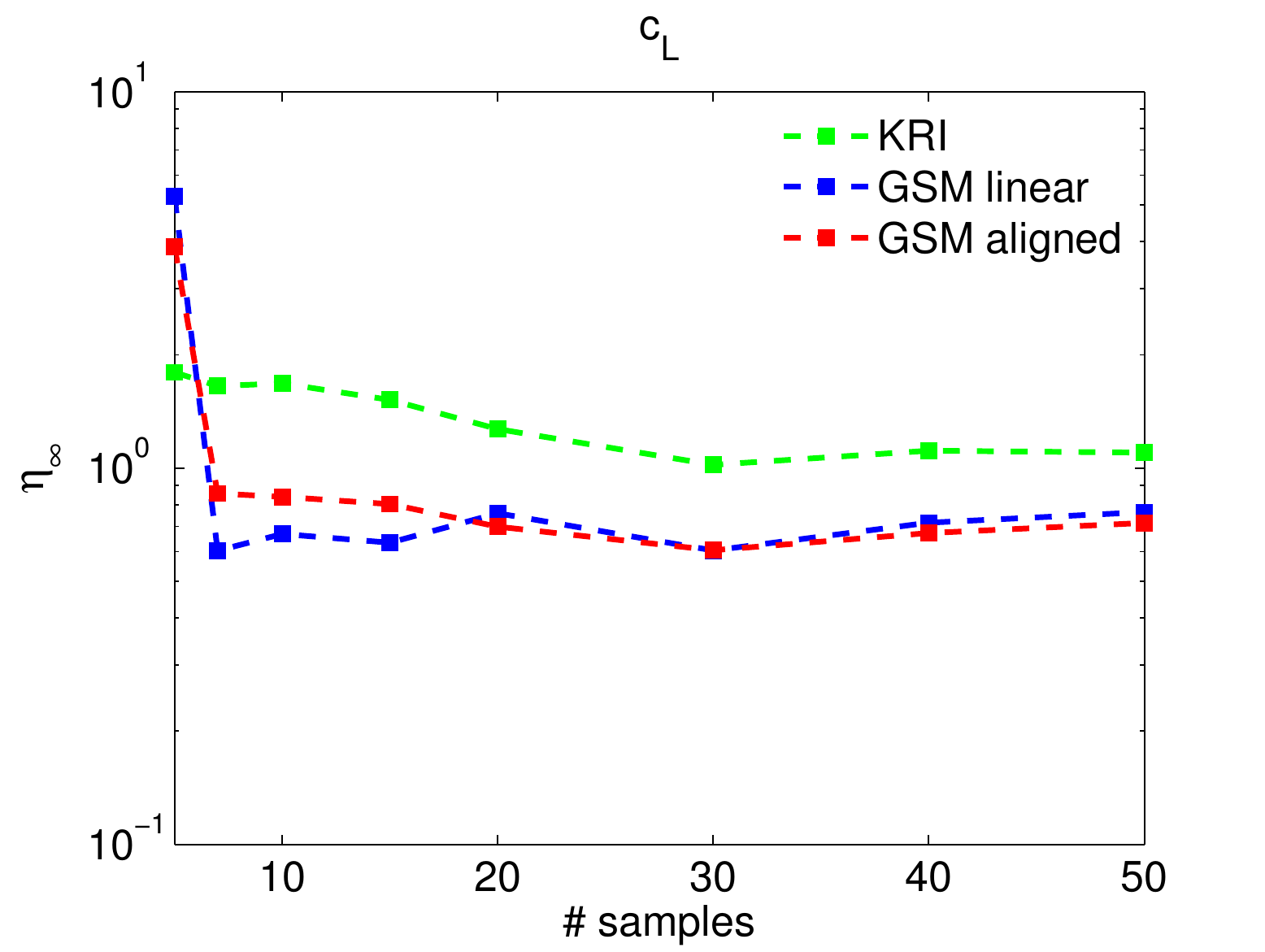}\\
  \includegraphics[width=0.49\textwidth]{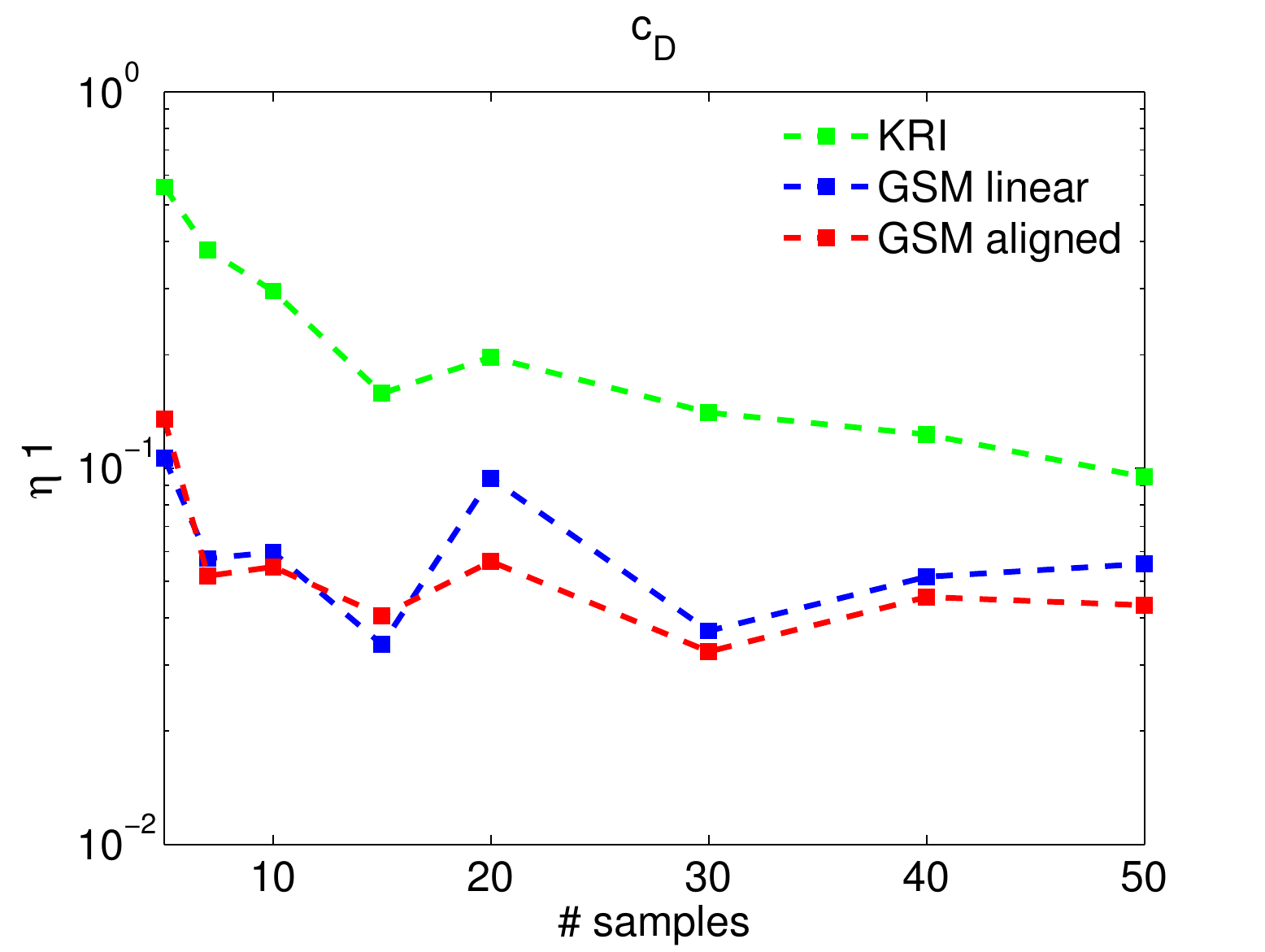} 
  \includegraphics[width=0.49\textwidth]{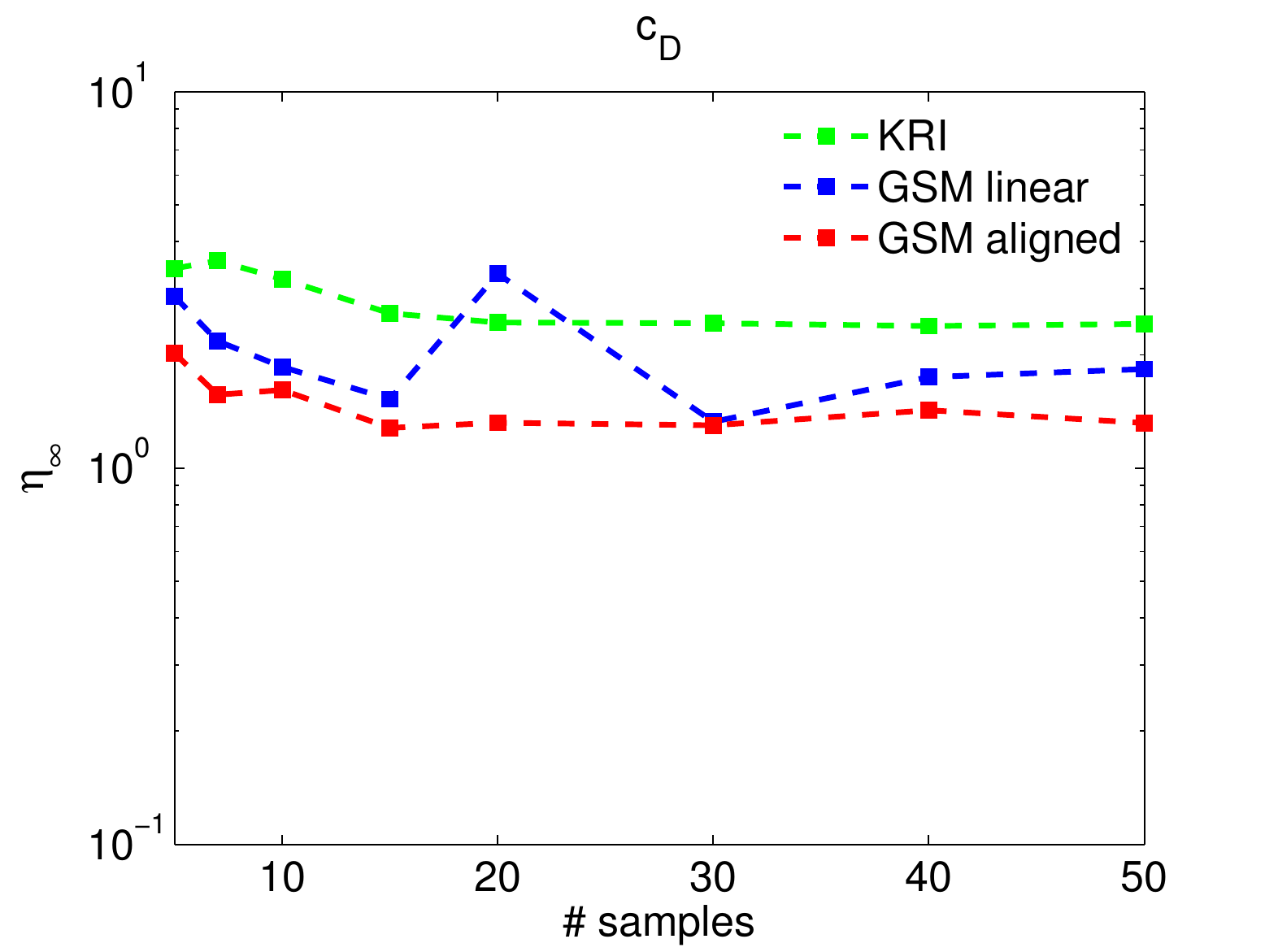}\\
  \includegraphics[width=0.49\textwidth]{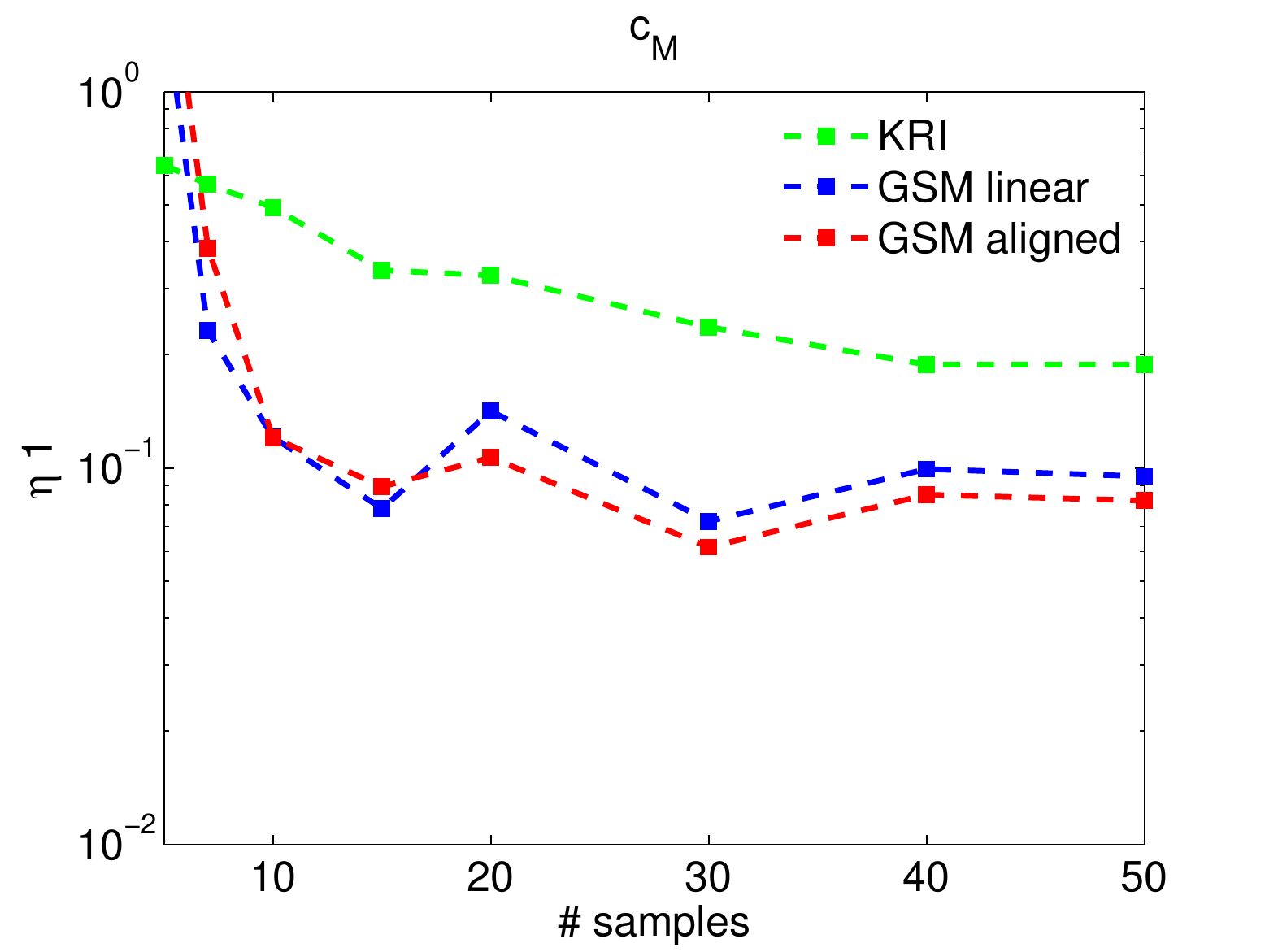} 
  \includegraphics[width=0.49\textwidth]{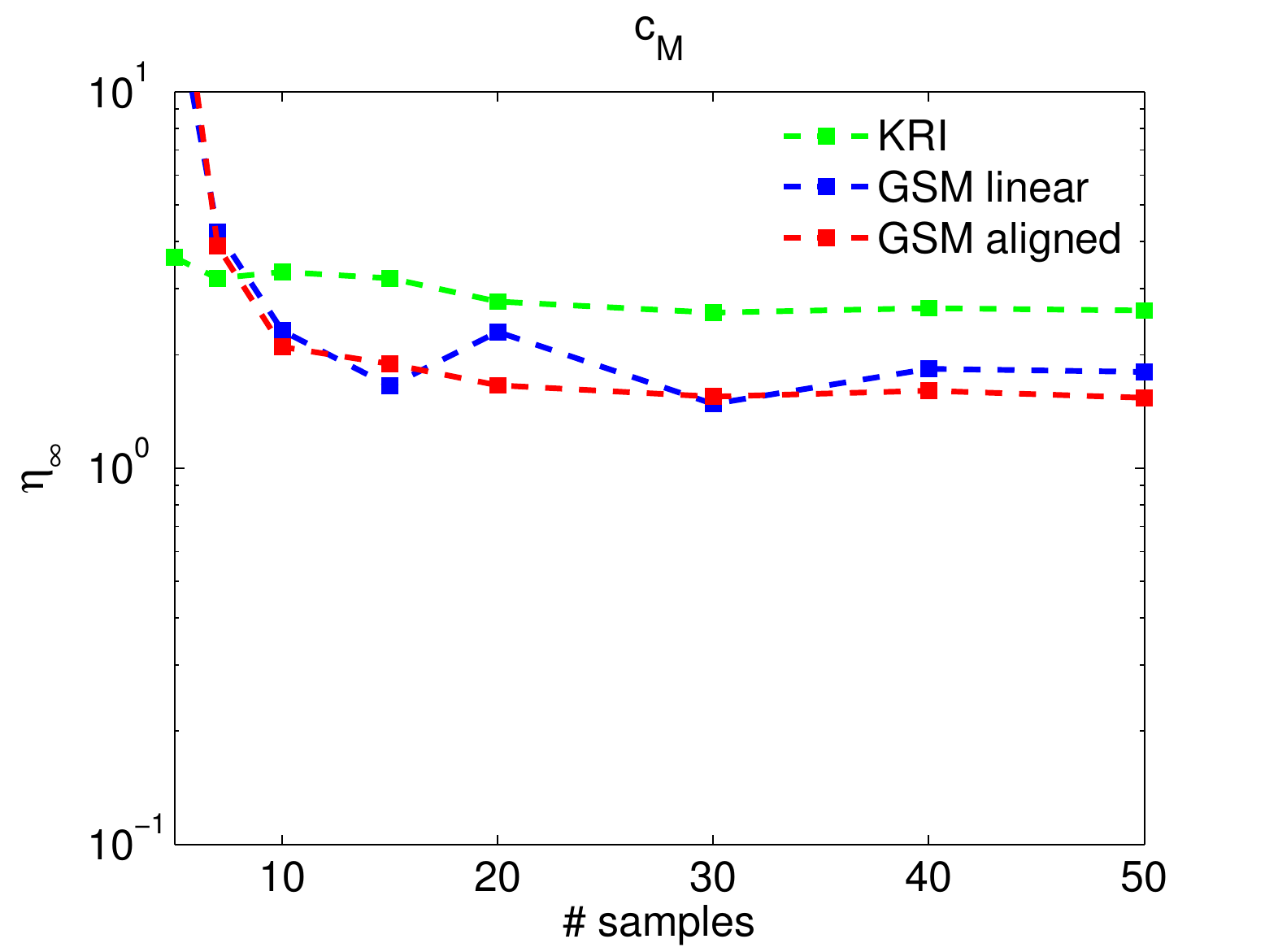} 
}
\caption{Average performances for Latin hypercube samplings. }
\label{fig:performanceLHC}
\end{figure}

With the POD-basis, we are now able to perform hierarchical Kriging interpolation (\ref{secHK:HK1}) based on the generic surrogate model (\ref{secGPOD:GPOD}) for new test cases, i.e.\ for responses of new airfoil geometries. We choose a NACA 3.375 2.875 1 2 profile, an airfoil which is not contained in the database, to demonstrate our surrogate modeling framework. The CFD solver is evaluated on a $40\times 40$ discrete grid $\Omega^{\text{val}}\subset\Omega$ to generate a set of validation data for the purpose of error evaluation. We use 
\begin{align}
  \eta_{1}&= \frac{1}{\sigma\left|\Omega^\text{val}\right|} \sum_{x\in\Omega^\text{val}}\left|\widehat\phi(x)-\phi(x) \right| \\
  \text{and}\quad  \eta_{\infty}&= \frac{1}{\sigma}\max_{x\in\Omega^\text{val}}\left|\widehat\phi(x)-\phi(x) \right|,
\end{align}
the relative average error and the relative maximum error as measures of accuracy, where $\sigma$ denotes the standard deviation of the set $\left\{\phi(x)\right\}_{x\in\Omega^\text{val}}$ and $\left|\Omega^\text{val}\right|=1600$. In a first study, the performance of ordinary Kriging interpolation is compared to hierarchical Kriging based on the generic surrogate model with and without the alignment (figure \ref{fig:performanceLHC}). We compute 10 Latin hypercube samplings each for the sample sizes 5, 7, 10, 15, 20, 30, 40 and 50 and use them to generate interpolations, for each of the three surrogate modeling techniques and for each of the three responses $c_{l}$, $c_{d}$, $c_{m}$. For each sample size, an average performance of the 10 Latin hypercube samplings is computed. The figure shows clearly that both hierarchical Kriging methods (blue and red) outperform the ordinary Kriging method (green) in terms of average error as well as maximum error. Exceptions are the accuracy of the $c_{l}$ and the $c_{m}$ approximations for the Latin hypercube sample sizes of 5 and 7, respectively. The reason is that for very small sample sizes, the gappy POD approximation (\ref{secGPOD:minsum})/(\ref{secGPOD:mintransform2}) has too many degrees of freedom (number of basis elements plus optionally number of transformation parameters) compared to the number of conditions (number of samples). This results in the few samples being approximated very closely by the linear combination of (aligned) POD-basis elements at the cost of an unfavorable behavior in the rest of the domain $\Omega$. One can further observe that using transformations in the generic surrogate model (red) decreases the approximation errors compared to hierarchical Kriging based on non aligned generic surrogate models (blue) for sample sizes $\geq 20$. So the effort for the solution of a nonlinear least squares gappy POD problem is justified. The computational costs of surrogate model generation are negligible compared to one single evaluation of the CFD solver: the maximum CPU time for a sample size of 50 was 40 seconds. Comparing the performances of both hierarchical Kriging methods depending on the sample size, one recognizes that with increasing number of samples the accuracy is not necessarily improved, especially for 40 and 50 samples.

\begin{figure}[h!]
\centering{ 
  \includegraphics[width=0.49\textwidth]{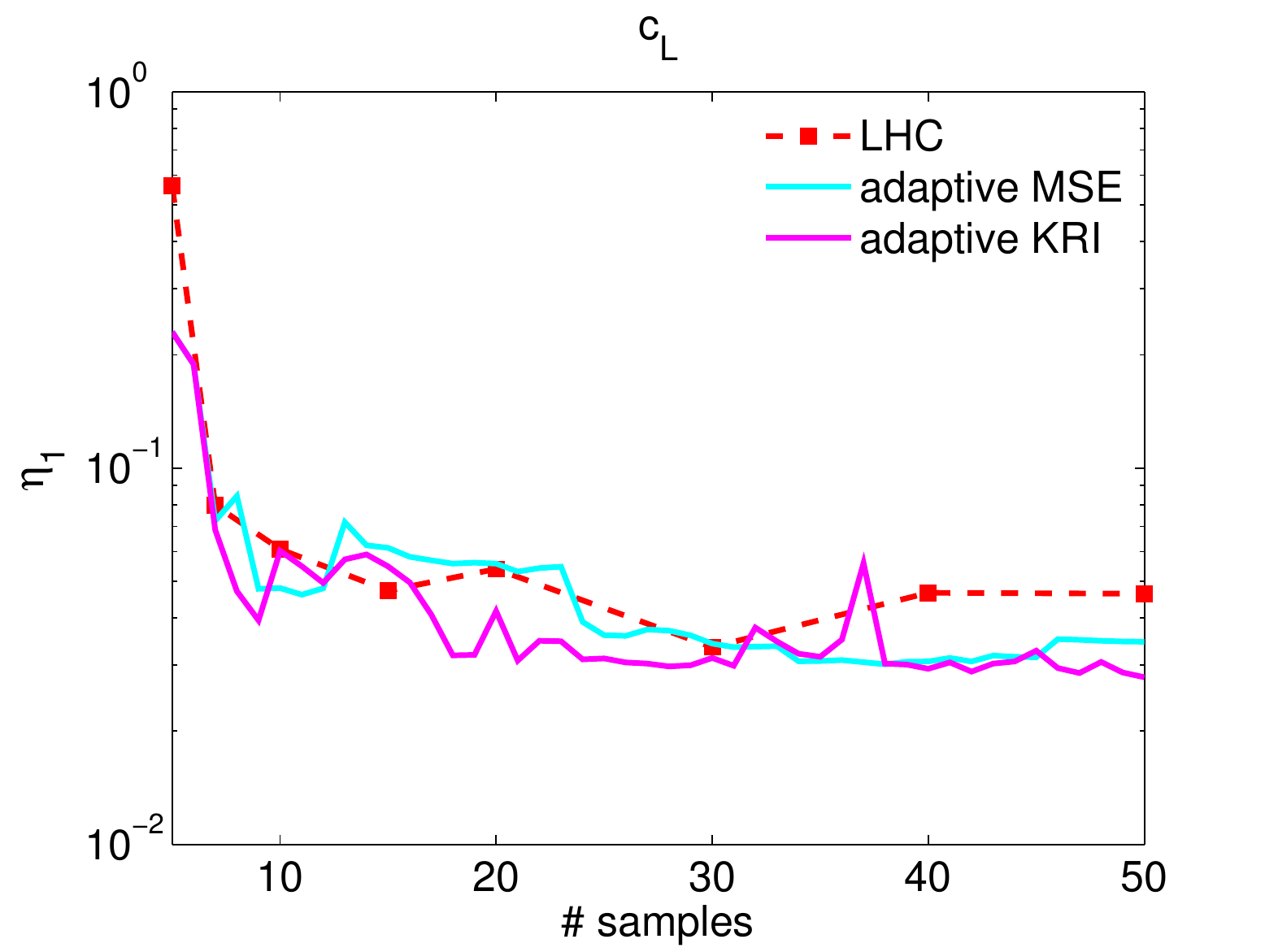} 
  \includegraphics[width=0.49\textwidth]{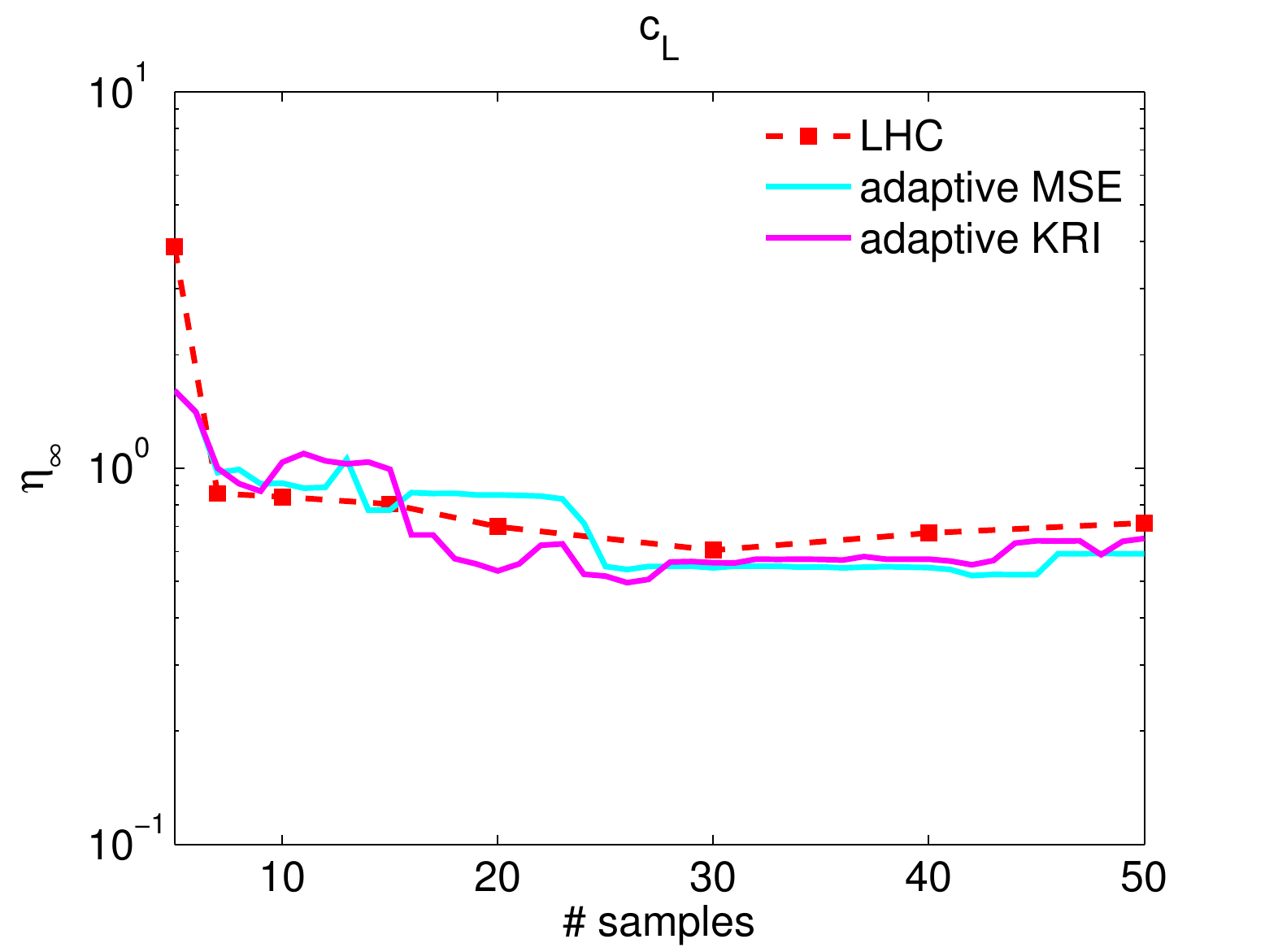}\\
  \includegraphics[width=0.49\textwidth]{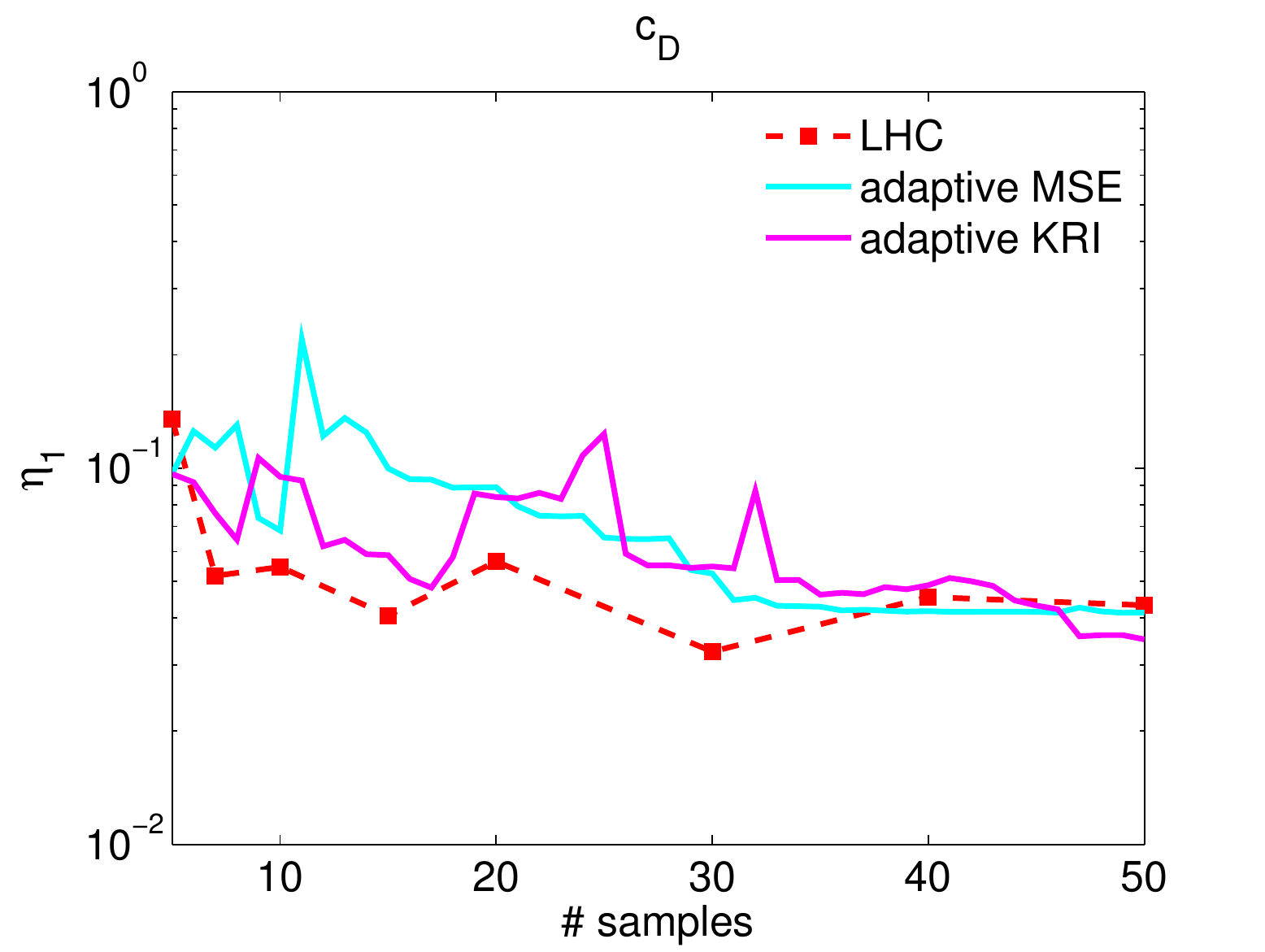} 
  \includegraphics[width=0.49\textwidth]{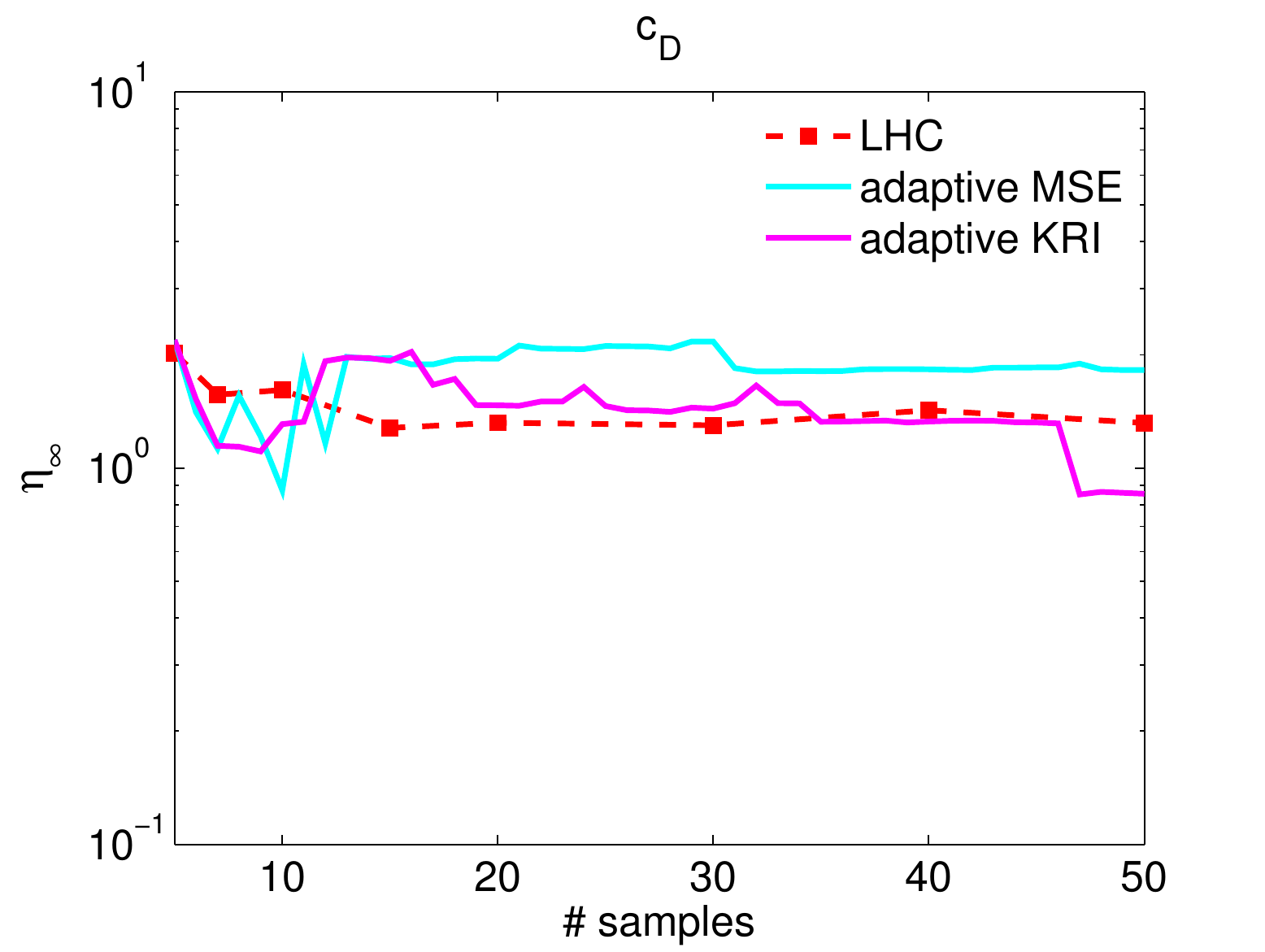}\\
  \includegraphics[width=0.49\textwidth]{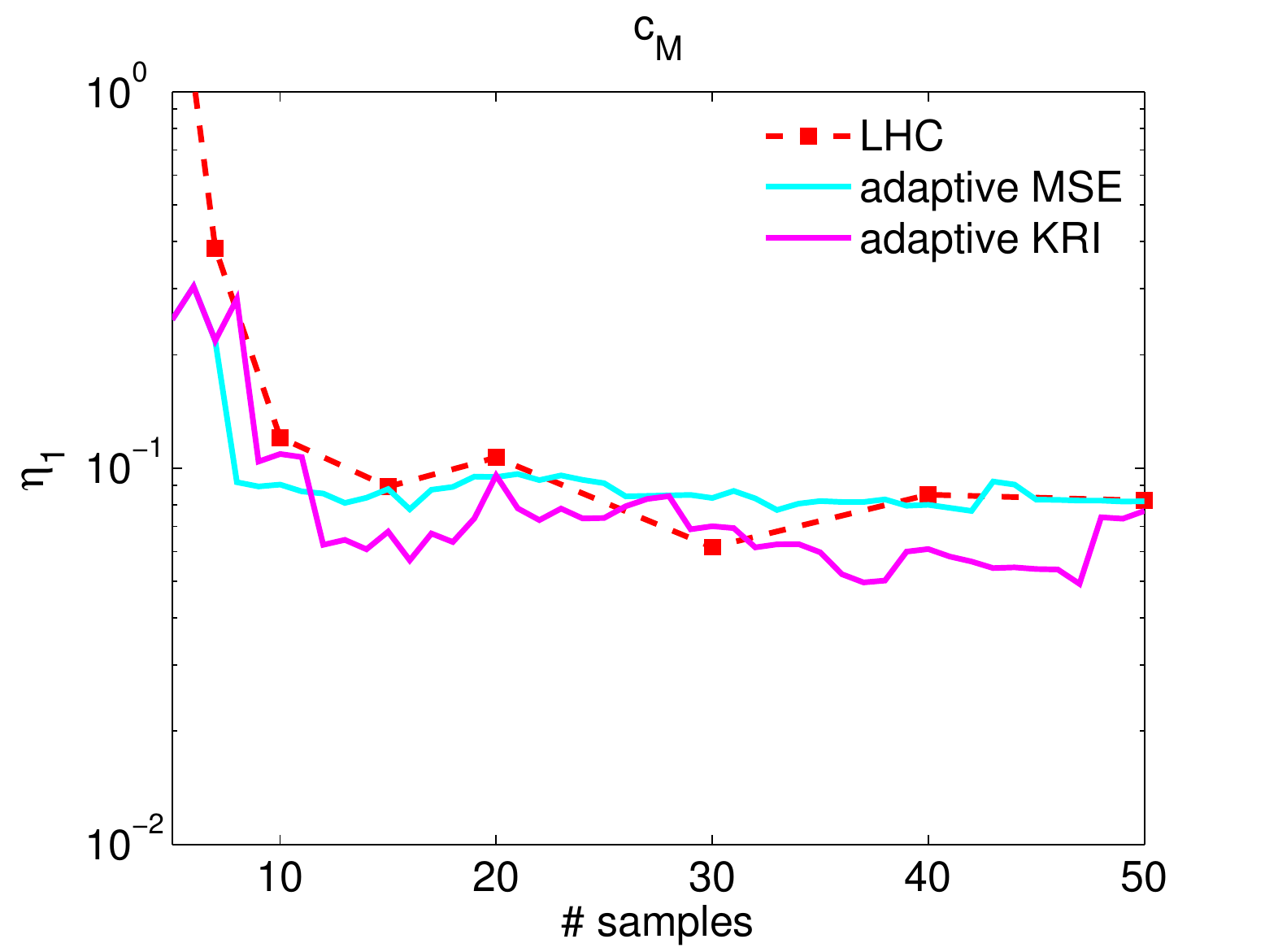} 
  \includegraphics[width=0.49\textwidth]{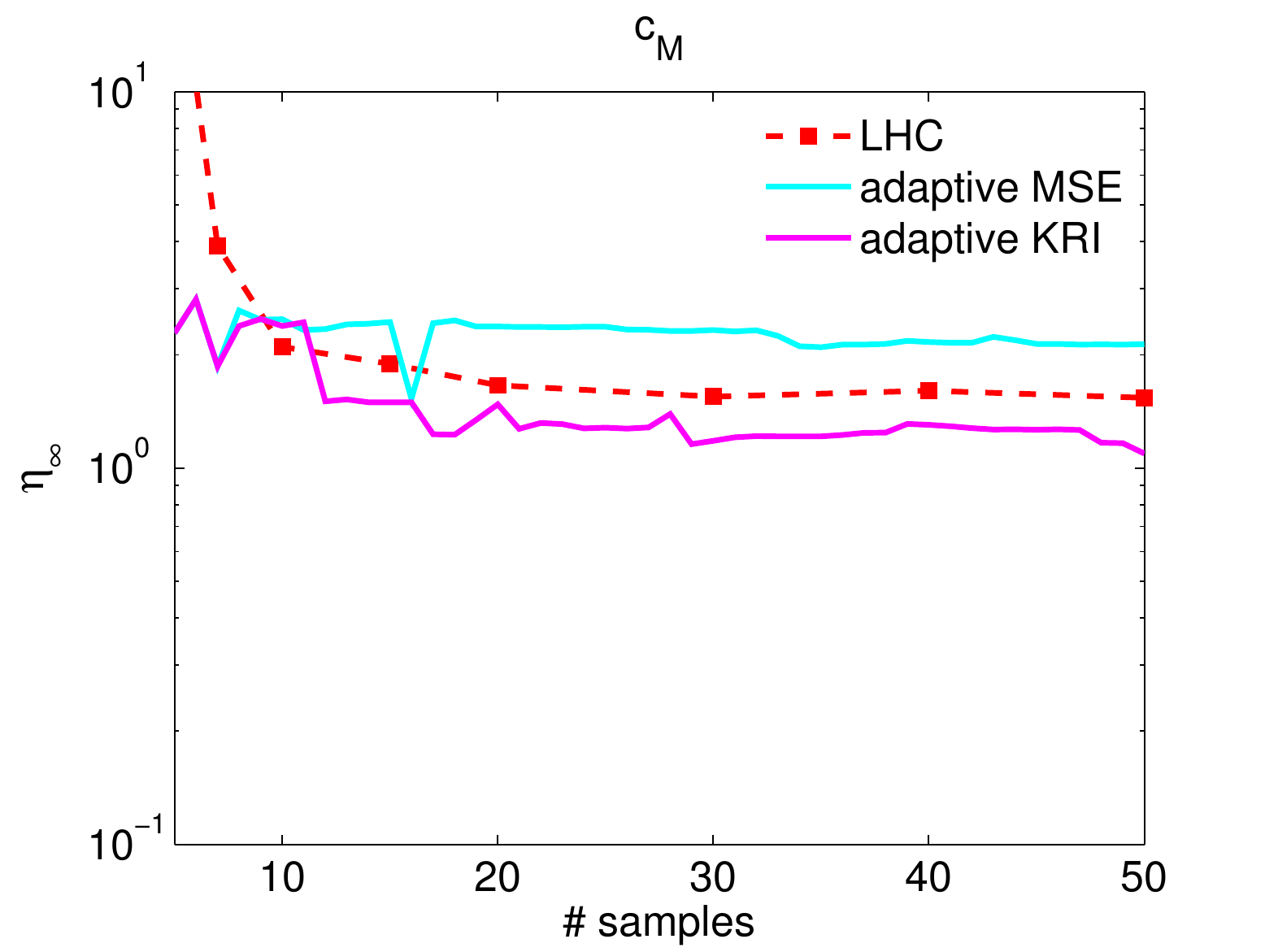} 
}
\caption{Comparison of Latin hypercube samplings and adaptive sampling strategies.}
\label{fig:performanceadaptive}
\end{figure}
\begin{figure}[h!]
\centering{ 
  \includegraphics[width=0.49\textwidth]{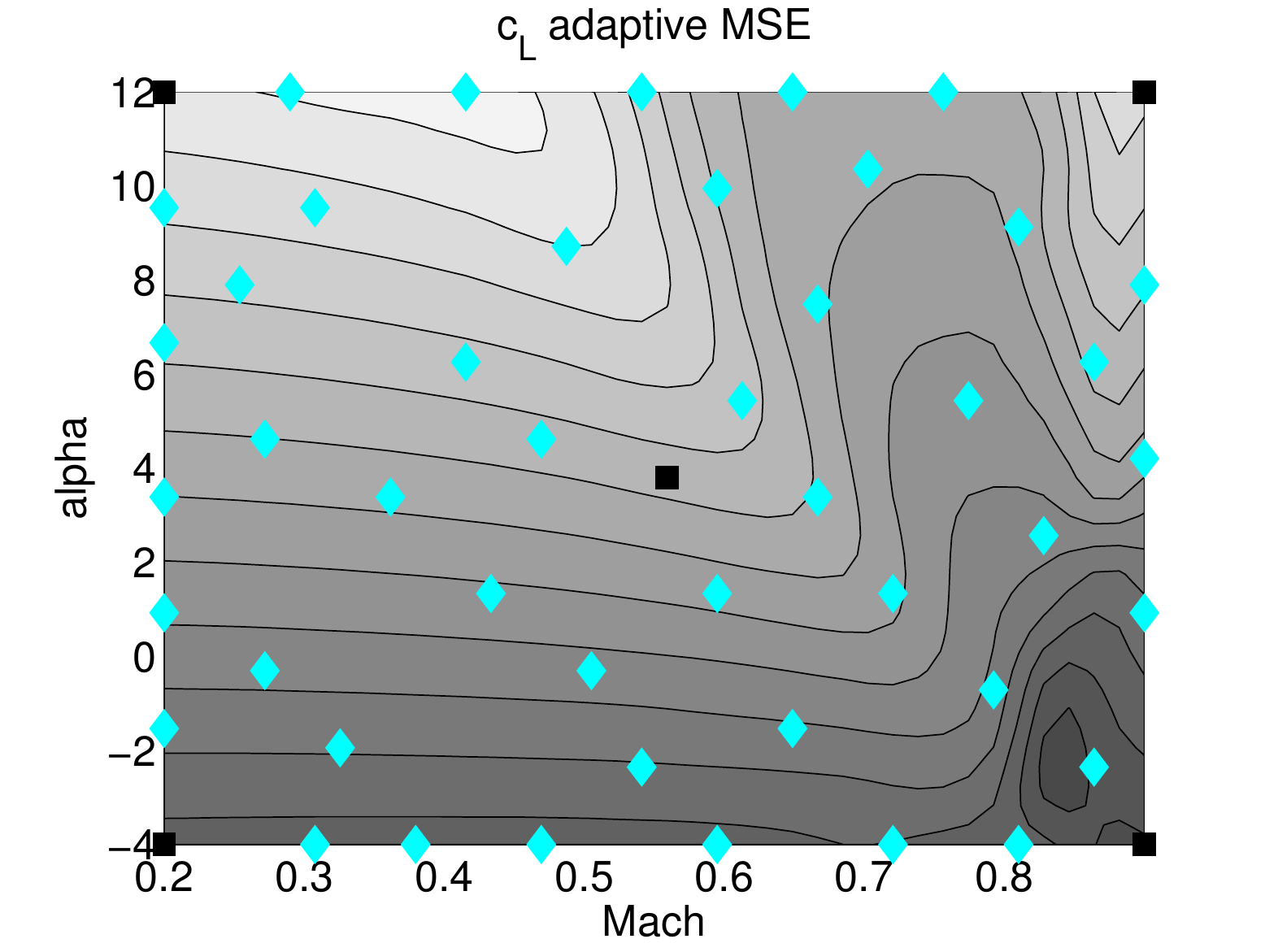} 
  \includegraphics[width=0.49\textwidth]{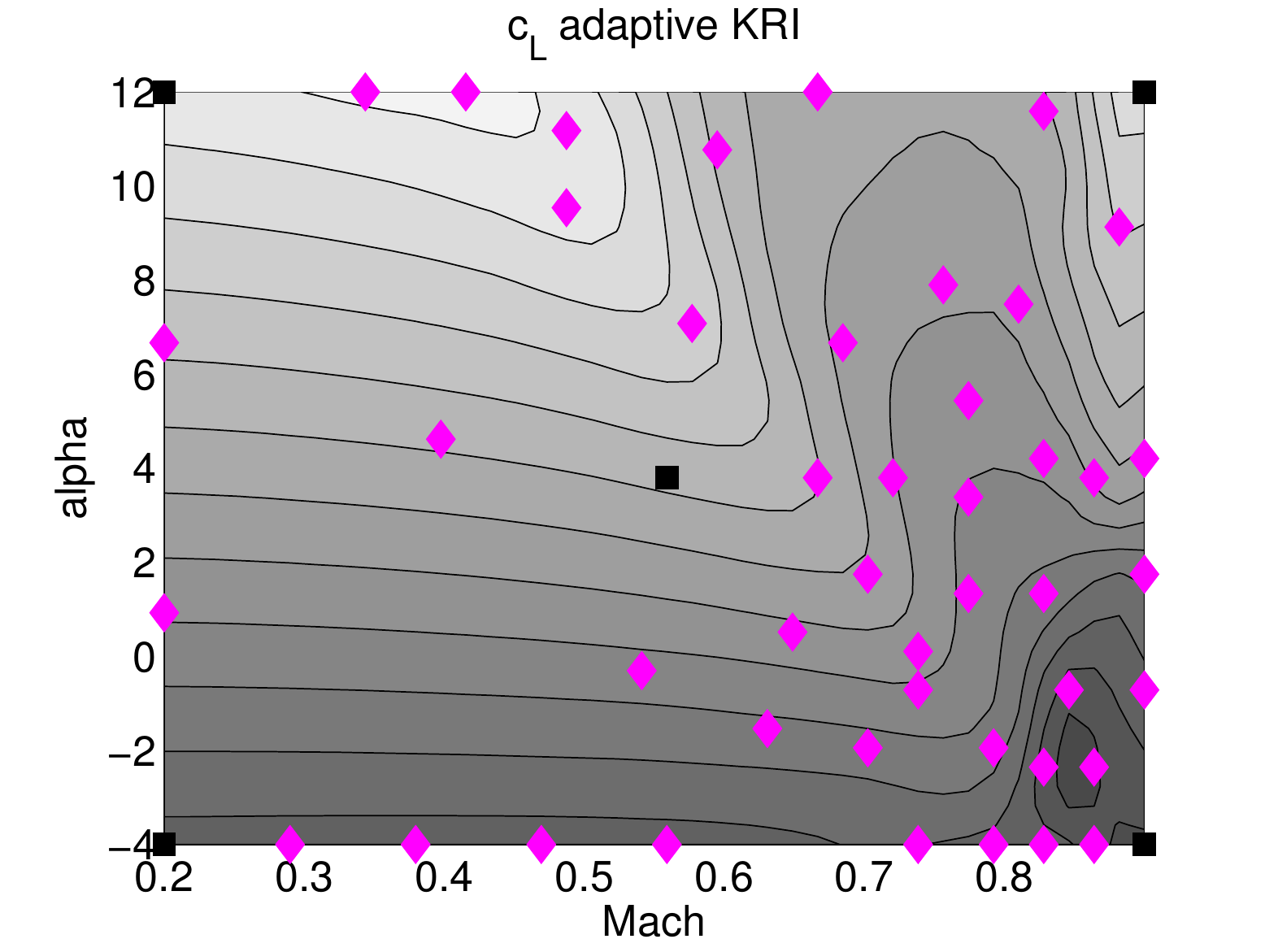}\\
  \includegraphics[width=0.49\textwidth]{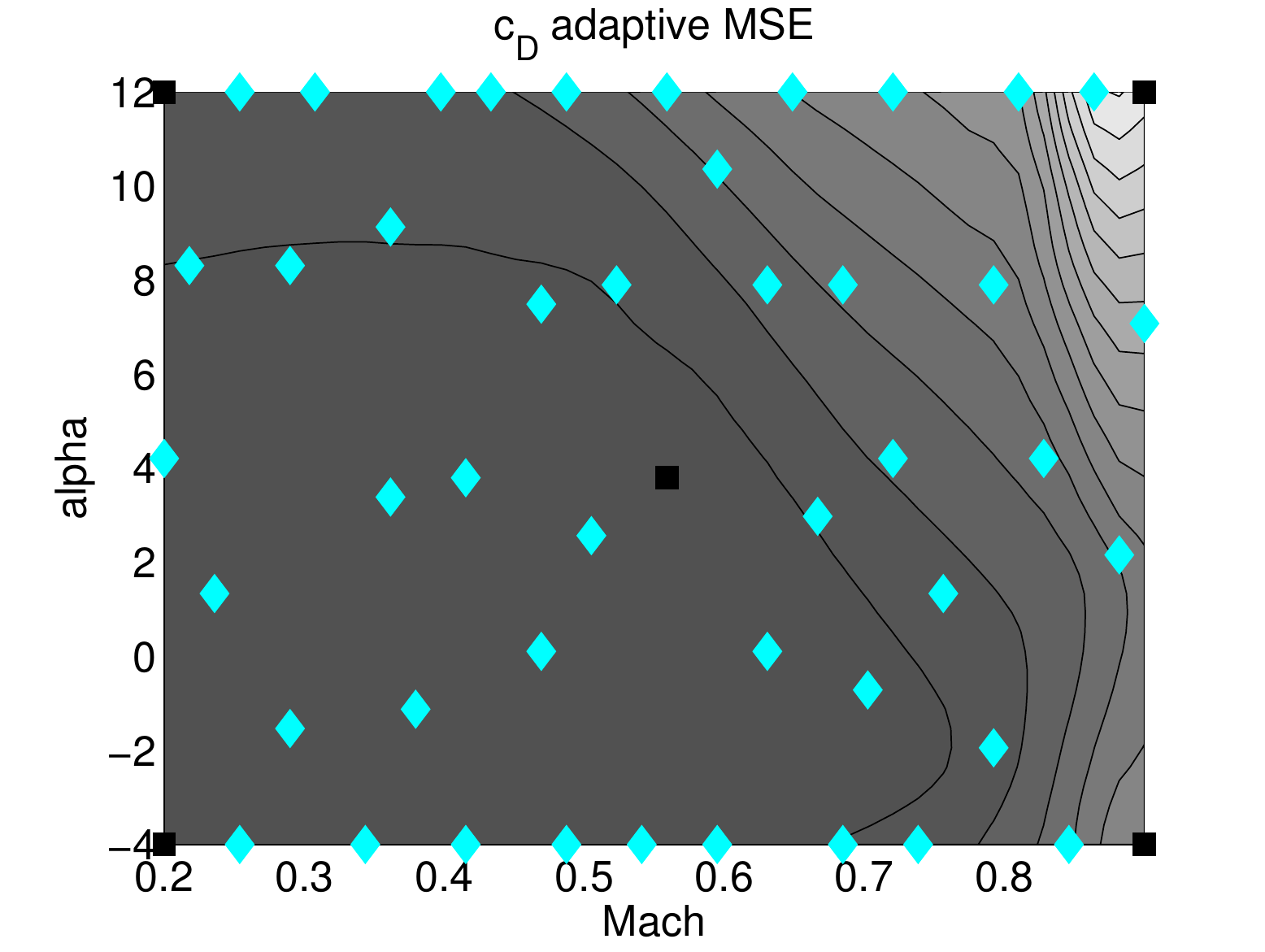} 
  \includegraphics[width=0.49\textwidth]{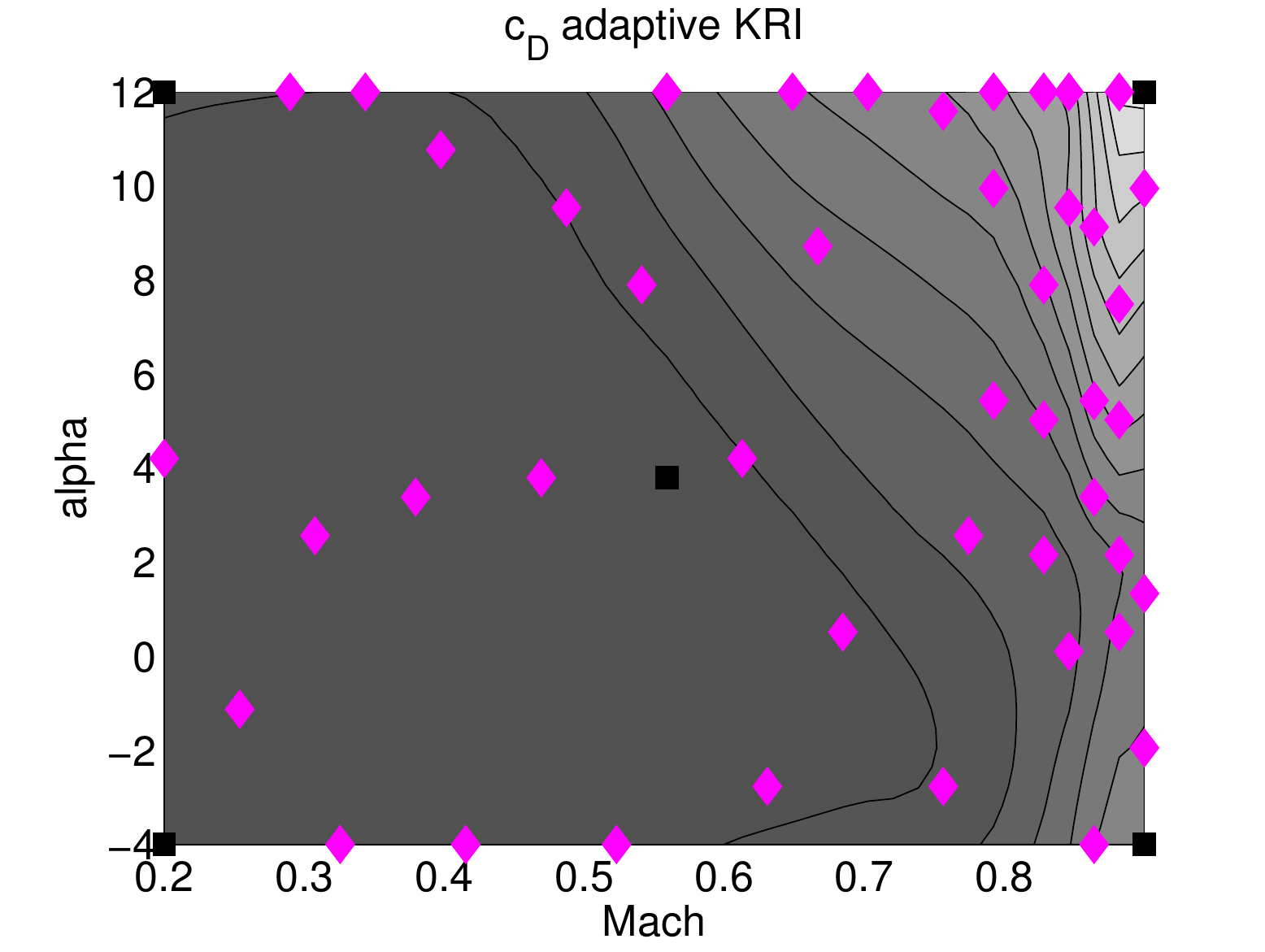}\\
  \includegraphics[width=0.49\textwidth]{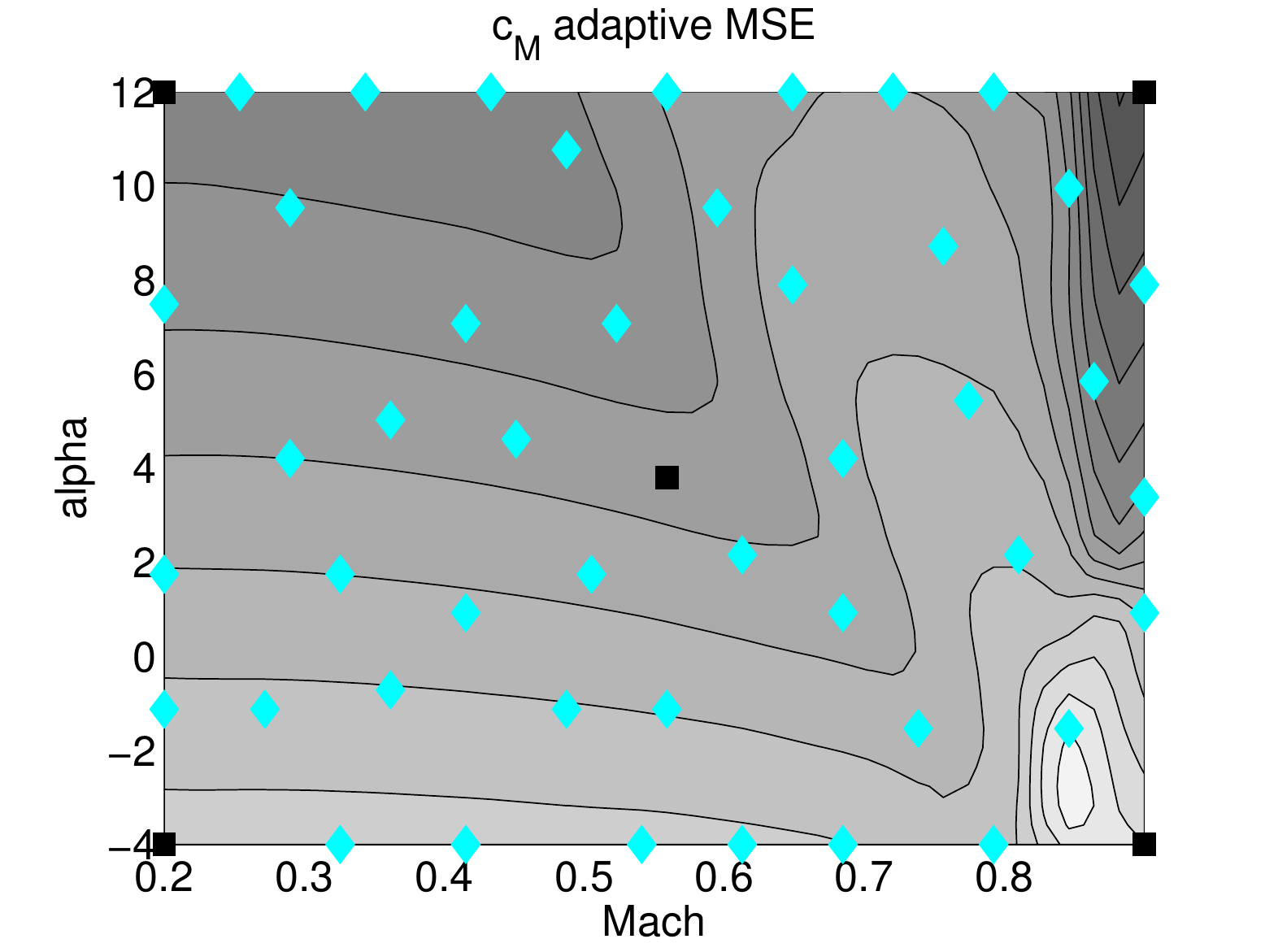} 
  \includegraphics[width=0.49\textwidth]{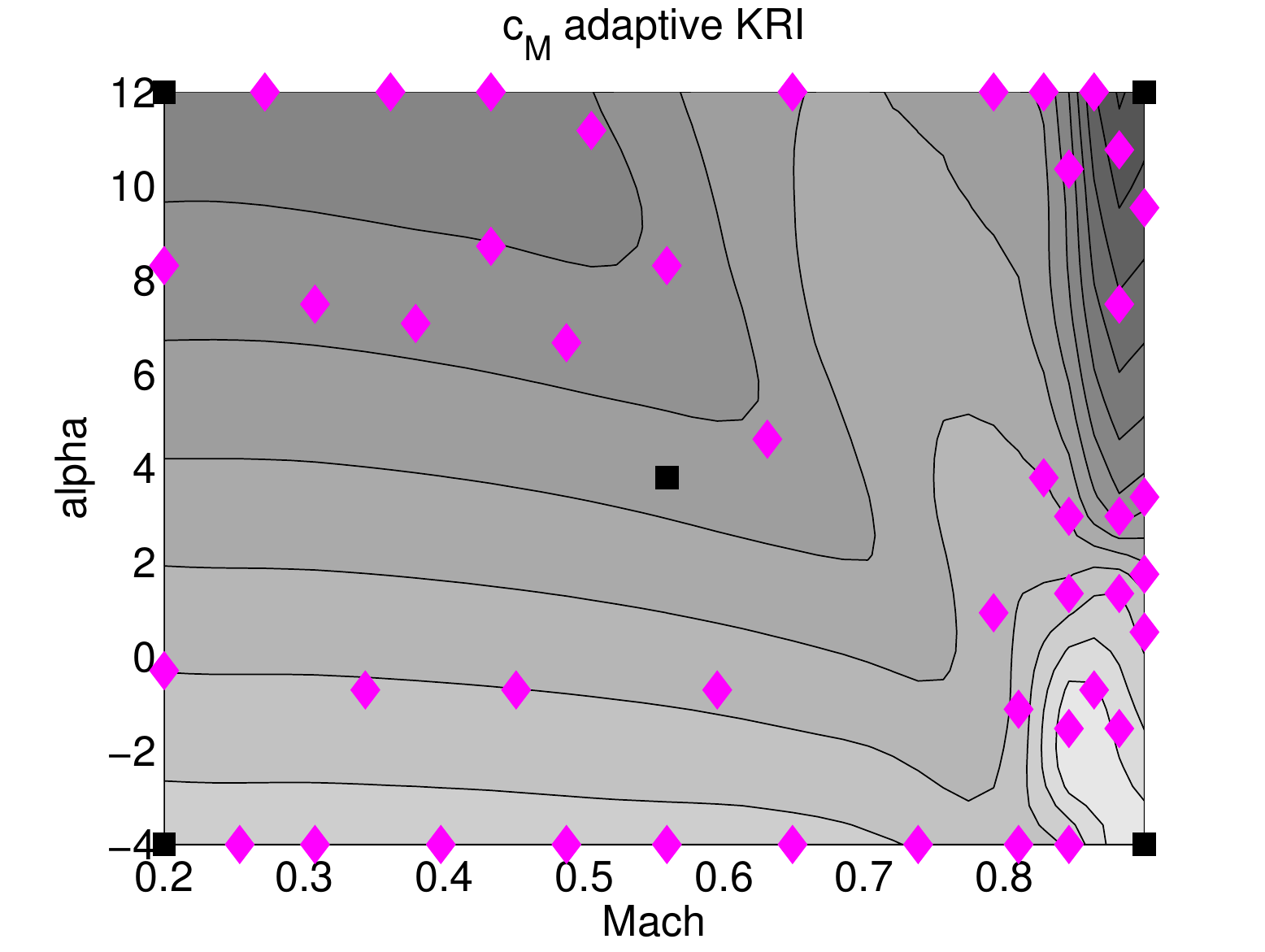} 
}
\caption{Distribution of samples.}
\label{fig:distributionadaptive}
\end{figure}

Therefore, we also apply adaptive sampling strategies. Unlike Latin hypercube or Monte Carlo sampling strategies, the samples are generated sequentially. At every stage of the adaptive process, a surrogate model is generated and assessed in order to find a new sample location $x^{*}$, in which the unknown function is evaluated. The data pair $(x^{*},\phi(x^{*}))$ is subsequently added to the current sampling. We investigated two adaptive sampling strategies, which are easy to implement in the generic surrogate modeling framework. The first one chooses the new sample $x^{*}$ where the predicted mean squared error $\MSE\left[\widehat\phi (x)\right]:= \E\left[\left(\widehat\phi (x)-\phi(x)\right)^{2}\right]$ is highest \cite{BR2011}. This quantity can easily be computed by
\begin{equation}
\MSE\left[\widehat\phi (x)\right] = \sigma^{2}\left(1- \begin{pmatrix} r(x) \\ \widetilde{\overline \phi}(\overline{\overline x}(p),p,a )\end{pmatrix}^{\top} \begin{bmatrix}R & \Phi \\\Phi\T & 0 \end{bmatrix}^{-1} \begin{pmatrix} r(x) \\ \widetilde{\overline \phi}(\overline{\overline x}(p),p,a )\end{pmatrix} \right),
\end{equation}
it is equal to zero in every previously sampled location and grows with distance to them.
The new sample location is determined by
\begin{equation}
  x^{*}=\operatorname{arg}\max_{x\in\Omega^{\text{val}}}\MSE\left[\widehat\phi (x)\right].
\end{equation}
A second adaptive sampling strategy for variable fidelity modeling is proposed in \cite{HanGH2010}. In the sequel, we distinguish between the GSM-based hierarchical Kriging $\widehat\phi_{\text{GSM}}(x):=\widehat\phi(x)$ and an ordinary Kriging interpolation $\widehat\phi_{\text{KRI}}(x)$. Assuming that both surrogate models will converge to the true function $\phi(x)$ with growing (maybe very large) number of samples, a new sample is added where the error between both models is highest:
\begin{equation}
  x^{*}=\operatorname{arg}\max_{x\in\Omega^{\text{val}}}\left|\widehat\phi_{\text{GSM}}(x) - \widehat\phi_{\text{KRI}}(x)\right|.
\end{equation}

For the hierarchical Kriging based on aligned gappy POD fittings which performed best in the study above, we compare both adaptive sampling strategies to the average Latin hypercube performance (red) in figure \ref{fig:performanceadaptive}. For both methods we use an initial sampling with 5 points. For the MSE method (cyan) no clear assessment is possible. For the $c_{l}$ response its performance is comparable to the Latin hypercube samplings up to 30 samples and outperforms them for larger sample sizes. The $c_{d}$ response is approximated more accurately by the Latin hypercube samplings than by the MSE method throughout almost all sample sizes, only for 40-50 samples the adaptive MSE is more accurate in terms of the average error $\eta_{1}$. For the $c_{m}$ response, its average error is comparable to the nonadaptive samplings' performance, while the maximum error can not be reduced in the adaptive process and is higher than for Latin hypercubes. The second adaptive strategy (magenta) seems to produce more robust results. Despite a long starting phase for the $c_{d}$ response, where it performs worse than the Latin hypercube samplings, at the latest for 45 samples it is more accurate, even earlier for $c_{l}$ and $c_{m}$. Figure \ref{fig:distributionadaptive} illustrates surrogate models and the distributions of sample points generated by the two adaptive methods. The 5 point initial sampling is depicted in black squares. The MSE method merely generates a space filling design. When the correlation lengths of the $x^{1}$ and the $x^{2}$ dimension differ, many samples are placed at the lower and upper boundary of the axis with the larger correlation length, an observation already described in a previous publication \cite{BR2011}. The second method adds more samples in regions where the response has larger gradients or curvature, e.g.\ high Mach number for all three responses or the already mentioned curvature of the $c_l$ response. This behavior is desirable, since traditionally in these regions surrogate models are least accurate and need more samples to describe the characteristics of the true function $\phi(x)$. On the other hand, the lower and upper boundary of the $\alpha$-axis obtain almost as many samples as in the MSE sampling strategy. These two attributes make this sampling strategy the most powerful in our test cases. We emphasize that in a previous study of similar test cases \cite{BR2011}, adaptive sampling strategies also required a starting phase and the benefit compared to Latin hypercube samplings was observed only after a total number of 40-60 samples.

\section{Conclusions}\label{sec:con}
In this paper a new approach for globally valid surrogate models was developed. We proposed a variable fidelity framework, which uses a generic surrogate model as a global trend. Generic surrogate models can be used, whenever surrogate models for multiple test cases of a predefined problem class are available. We addressed how to establish correspondence between these database functions and how to identify characteristic structures by POD. The generic surrogate model was introduced as a gappy POD fit to the sample data and we extended the method to nonlinear transformations. Hierarchical Kriging, a recently developed VFM method, was used to interpolate sample data based on the generic surrogate model. In contrast to other VFM methods, where high-fidelity and low-fidelity data are both computed ``online'', i.e.\ for every new test case of a problem class, in generic surrogate modeling the database functions are computed once and stored (``offline''), such that for every new test case of the mutual problem class only the high-fidelity data has to be computed online. We validated the methods in a test case of aerodynamic coefficients of two-dimensional airfoils depending on the input parameters $(\textit{Ma},\alpha)$. A database of surrogates was generated using RANS computations with DLR TAU for 24 airfoils. For a new airfoil not contained in the database interpolations were generated for Latin hypercube samplings. Hierarchical Kriging based on generic surrogate models performed more accurate than ordinary Kriging interpolations and the benefit was largest for sample sizes up to 30. We also showed that further improvement of the approximation quality is possible using adaptive sampling strategies. Requiring significantly less expensive samples to achieve a desired accuracy than ordinary Kriging, hierarchical Kriging based on generic surrogate models describes an efficient framework in surrogate modeling.

\addcontentsline{toc}{chapter}{References} 


\end{document}